%% file: SOSC.tex
\documentclass[11pt]{amsart}
\usepackage{geometry}                
\usepackage{graphicx}
\usepackage{amssymb}
\usepackage{epstopdf}
\usepackage{hyperref}
\hypersetup{
    bookmarks=true,         
    unicode=false,          
    pdftoolbar=true,        
    pdfmenubar=true,        
    pdffitwindow=false,     
    pdfstartview={FitH},    
    pdftitle={My title},    
    pdfauthor={Author},     
    pdfsubject={Subject},   
    pdfcreator={Creator},   
    pdfproducer={Producer}, 
    pdfkeywords={keywords}, 
    pdfnewwindow=true,      
    colorlinks=true,       
    linkcolor=red,          
    citecolor=magenta,        
    filecolor=cyan,      
    urlcolor=blue           
}

\input{customcommands}

\newtheorem{lemma}{Lemma}
\newtheorem{corollary}{Corollary}
\newtheorem{theorem}{Theorem}
\newtheorem{example}{Example}
\newtheorem{algorithm}{Algorithm}
\newtheorem{test}{SOSC Test}

\title{Hessian-Free Methods for Checking the Second-Order Sufficient Conditions in Equality-Constrained Optimization and Equilibrium Problems}
\author{W. Ross Morrow}
\address{Departments of Mechanical Engineering and Economics \\
	Iowa State University \\
	Ames, Iowa 50011}
\email{wrmorrow@iastate.edu}

\date{\today}                                           

\begin{document}
\maketitle


\begin{abstract}
Verifying the Second-Order Sufficient Condition (SOSC), thus ensuring a stationary point locally minimizes a given objective function (subject to certain constraints), is an essential component of non-convex computational optimization and equilibrium programming. This article proposes three new ``Hessian-free'' tests of the SOSC that can be implemented efficiently with gradient evaluations alone and reveal feasible directions of negative curvature when the SOSC fails. The Bordered Hessian Test and a Matrix Inertia test, two classical tests of the SOSC, require explicit knowledge of the Hessian of the Lagrangian and do not reveal feasible directions of negative curvature should the SOSC fail. Computational comparisons of the new methods with classical tests demonstrate the relative efficiency of these new algorithms and the need for careful study of false negatives resulting from accumulation of round-off errors.
\end{abstract}




\section{INTRODUCTION}

Verifying the Second-Order Sufficient Condition (SOSC) to certify local optimality of points computed by optimization software is an important but underdeveloped component of computational optimization. Existing optimization solvers compute points satisfying a First-Order or Second-Order Necessary Condition (FONC/SONC), typically without checking the corresponding SOSC. As a result, there remain cases in which computed points are {\em not} optimizers such as, for example, Example \ref{EX:Cube}. 

While verifying the SOSC is important for non-convex optimization, verifying the SOSC is {\em essential} for equilibrium models increasingly being employed in Economics, Operations Research, and Engineering. While much of the theory of computing equilibria relies on assumptions of convexity to ensure the first-order conditions imply optimality rather than just stationarity \cite{Facchinei09}, examples of non-convex games are rapidly appearing in important applications. So long as equilibrium programming methods for general, non-convex games are restricted to solving the combined first-order conditions, algorithms can compute simultaneous first-order points but cannot distinguish equilibria from other types of first-order points; see, e.g., Example \ref{EX:EX1}. Verifying the SOSC is thus fundamental to properly computing equilibria. 

At least two ways to test the SOSC have been known for some time. The classical ``Bordered Hessian Test'' (BHT) \cite{Gillespie51, McFadden78, Papalambros00, Luenberger09} can, in principle, be used to verify or reject the SOSC at points computed by optimization software. Computationally implementing the BHT relies on a set of nested LU factorizations that can be efficiently taken with LU factorization updating. Determining the inertia of a ``KKT matrix'' \cite{Gould85, Nocedal06} is another ``classical'' way to verify or reject the SOSC. Efficient implementations of this ``Inertia test'' are available in the form of existing factorization packages that can compute matrix inertia, and may be easily integrated into SQP solvers for constrained optimization problems. 

Both of these ``classical'' tests require an explicit representation of the Hessian of the Lagrangian, which may not be available. Popular optimization solvers including \texttt{MINOS}, \texttt{LANCELOT}, \texttt{SNOPT}, \texttt{KNITRO}, and \texttt{matlab} can utilize gradient evaluations alone through quasi-Newton updates, finite-differences, or automatic differentiation \cite{Nocedal06}. The Hessian of the Lagrangian will thus not be available for verifying the SOSC when gradient evaluations alone are used. ``Hessian-free'' algorithms requiring only Hessian-vector products \cite{Nocedal06} could seamlessly integrate SOSC checks with optimization and equilibrium solvers that do not require users to provide formulas for the second derivatives of the objective or the constraint functions. Such algorithms can fully exploit sparsity patterns in the objective and constraint functions, as well as be implemented using directional finite differences to approximate Hessian-vector products. 

This article presents three new Hessian-free algorithms for verifying or rejecting the SOSC at first-order points in smooth equality-constrained optimization or equilibrium problems. The first algorithm is based on Cholesky factorization, the most efficient and stable method for testing Hessian positive-definiteness without constraints \cite{Golub96, Higham02}. The second algorithm is based on an ``oblique'' Gram-Schmidt orthogonalization, and the third algorithm is a modification of the Projected Conjugate Gradient algorithm developed for constrained quadratic programming \cite{Gould01, Nocedal06}. By being Hessian-free, these algorithms can take full advantage of sparsity patterns in the Hessian of the Lagrangian when it is known and may significantly reduce the number of gradient evaluations necessary to verify or reject the SOSC when the Hessian is not available at all. 

Another important feature of an algorithm to verify the SOSC is how easily a feasible direction of negative curvature can be computed should the SOSC fail. Certain optimization algorithms use feasible directions of negative curvature to promote global convergence to second-order necessary points; see, e.g., \cite{Conn00}. Recovery of a feasible direction of negative curvature when an the SOSC fails enables ``warm restarts'' of optimization algorithms that may not already take advantage of second-order information. The classical BHT and Inertia tests do not provide an obvious path to computing such directions when the SOSC fails. In contrast, the Hessian-free algorithms proposed here make computation of such directions straightforward when the SOSC fails. Indeed, two of the Hessian-free algorithms presented here fail precisely by finding such directions. 

Finally, useful algorithms for verifying the SOSC should not return false-positives or false-negatives due to the accumulation of round-off errors. Unfortunately the numerical accuracy of verifying the SOSC with {\em any} approach has not yet been addressed. Section \ref{SEC:Examples} provides computational comparisons of the different algorithms including an example that illustrates significant potential for erroneous results due to round-off errors when the constraint gradients are very nearly linearly dependent. In particular, no method can be considered accurate for certain problems with very nearly linearly independent constraint gradients, even with relatively few variables and constraints. 


\section{OPTIMIZATION AND EQUILIBRIUM PROBLEMS}

\subsection{Optimization Problems}

This article considers the equality-constrained, continuous variable optimization problem
\begin{equation}
	\label{EQN:Opt}
	\begin{aligned}
		\mnm 	&\quad f(\vec{x}) \\
		\wrt 		&\quad \vec{x} \in \R^N \\
		\sto 		&\quad \vec{c}(\vec{x}) = \vec{0}
	\end{aligned}
\end{equation}
where $f: \R^N \to \R$, $\vec{c} : \R^N \to \R^M$, and $N,M \in \N, \; M < N$. The objective ($f$) and constraints ($\vec{c}$) are nonlinear and twice continuously differentiable functions of $N$ variables $\vec{x} = (x_1,\dotsc,x_N)$. The component functions of $\vec{c}$ are denoted by $c_i$, for $i \in \{1,\dotsc,M\}$. 

\begin{theorem}(Optimality Conditions \cite[Theorems 12.1, 12.5, 12.6]{Nocedal06})
	\label{THM:OptCond}
	Consider Problem (\ref{EQN:Opt}), and assume the Linear Independence Constraint Qualification (LICQ) is satisfied \cite[Def. 12.4]{Nocedal06}. Define the ``Lagrangian'' $\set{L}(\vec{x},\bsym{\lambda}) = f(\vec{x}) - \bsym{\lambda}^\top\vec{c}(\vec{x})$, and let $\vec{A}(\vec{x}) = D\vec{c}(\vec{x}) \in \R^{M \times N}$ denote the Jacobian matrix of the constraint function evaluated at $\vec{x} \in \R^N$. 
	\newline{\em {\bfseries FONC:}} 
	If $\vec{x}^* \in \R^N$ is a local solution to Problem (\ref{EQN:Opt}), then there exists some $\bsym{\lambda}^* \in \R^M$ such that $\nabla\set{L}(\vec{x}^*,\bsym{\lambda}^*) = \vec{0}$. That is, $\nabla f(\vec{x}^*) = \sum_{i=1}^M \lambda_i^* \nabla c_i(\vec{x}^*)$ and $\vec{c}(\vec{x}^*) = \vec{0}$. 
	\newline{\em {\bfseries SONC:}} 
	Moreover, $\vec{w}^\top\;\vec{H}(\vec{x}^*,\bsym{\lambda}^*)\;\vec{w} \geq 0$ for all $\vec{w} \in \set{C}(\vec{x}^*) = \{ \; \vec{h} : \vec{A}(\vec{x}^*)\vec{h} = \vec{0} \; \}$ where $\vec{H}(\vec{x},\bsym{\lambda})$ denotes the Hessian of the Lagrangian (or simply ``Hessian''):
	\begin{equation*}
		\vec{H}(\vec{x},\bsym{\lambda})
			= \nabla_{x,x}^2\set{L}(\vec{x},\bsym{\lambda})
			= \nabla_{x,x}^2f(\vec{x}) 
				- \sum_{i=1}^M \lambda_m \nabla_{x,x}^2c_i(\vec{x}). 
	\end{equation*}
	\newline{\em {\bfseries SOSC:}} 
	On the other hand, suppose $\vec{x}^* \in \R^N$ and $\bsym{\lambda}^* \in \R^M$ satisfy the FONC and $\vec{w}^\top\;\vec{H}(\vec{x}^*,\bsym{\lambda}^*)\;\vec{w} > 0$ for all $\vec{w} \in \set{C}(\vec{x}^*)$. Then $\vec{x}^*$ is an isolated local solution to Problem (\ref{EQN:Opt}).
\end{theorem}
Note that this SOSC also applies to {\em inequality} constrained optimization problems at {\em strictly complementary} stationary points \cite[Def. 12.5]{Nocedal06}, and as the continuous part of the optimality conditions for {\em mixed-integer} nonlinear optimization problems. The remainder of this article assumes that $\vec{x}^*$ and $\bsym{\lambda}^*$ satisfy the FONC and denotes $\vec{H}(\vec{x}^*,\bsym{\lambda}^*)$, $\vec{A}(\vec{x}^*)$, and $\set{C}(\vec{x}^*)$ by simply $\vec{H}$, $\vec{A}$, and $\set{C}$, respectively. Below the SOSC is denoted compactly using the symbol $\vec{H} \succ_{\set{C}} \vec{0}$. 

Many existing codes for solving problem (\ref{EQN:Opt}) solve a variant of the FONC without verifying the SOSC at computed points \cite{Nash98}. {\em Sequential Quadratic Programming} (SQP) methods \cite[Chapter 18]{Nocedal06} solve a sequence of local quadratic model problems that, in the equality-constrained case, corresponds to applying Newton's method to solve the FONC. SQP methods are currently implemented in \texttt{NPSOL} \cite{Gill01}, \texttt{SNOPT} \cite{Gill05}, \texttt{filterSQP} \cite{Fletcher02a, Fletcher02b, Fletcher02c}, \texttt{KNITRO} \cite{Byrd04,Byrd06}, and \texttt{matlab} \cite{Han77, Powell78a, Powell78b, Gill81, Fletcher87}. {\em Augmented Lagrangian} methods \cite[Chapter 17]{Nocedal06}, \cite{Bertsekas99} ``penalize'' the Lagrangian with a measure of the constraint violation and solve the FONC for a sequence of penalized problems with a variant of Newton's method. This method is currently implemented in the \texttt{MINOS} \cite{Murtagh78, Murtagh82, Murtagh98} and \texttt{LANCELOT} \cite{Conn92, LANCELOT} codes. Like Augmented Lagrangian methods, {\em Interior-Point} methods \cite[Chapter 19]{Nocedal06} for inequality-constrained problems solve the FONC for a sequence of approximate problems using a variant of Newton's method. \texttt{KNITRO} \cite{Byrd06}, \texttt{LOQO} \cite{Vanderbei99,Vanderbei06}, \texttt{IPOPT} \cite{Wachter06}, and \texttt{matlab} \cite{Byrd00,Byrd99} currently contain implementations of interior-point methods. 

Obtaining a solution to the FONC is {\em not} sufficient to declare the computed point a local solution to (\ref{EQN:Opt}), as shown in Example \ref{EX:Cube} below. In practice, sufficient decrease conditions on a merit function (or a filter mechanism) bias existing solvers towards computing constrained minimizers of $f$ \cite{Dennis96, Nocedal06, Schenk08}. Indeed, this bias towards optimizers is certainly one reason separate codes for verifying the SOSC do not currently exist. Sufficient decrease conditions certainly rule out some types of convergence: for example, these conditions rule out converge to local constrained {\em maximizers} of $f$. However, algorithms with sufficient decrease conditions can still converge to saddle points, as in Example \ref{EX:Cube}. 

\begin{example}
	\label{EX:Cube}
	Consider minimizing $f(x) = x^3$ over $\R$. $x^* = 0$ satisfies the FONC ($f^\prime(0) = 0$) and the SONC ($f^{\prime\prime}(0) = 0$), but does not satisfy the SOSC. Indeed, $f$ is not locally minimized at $x^* = 0$, or at any other finite $x$. Applying SQP \cite[Chapter 18.1]{Nocedal06} to this problem starting at any $x_0 > 0$ generates the sequence $x_n = x_{n-1}/2 = x_0/2^n \to 0$ as $n \to \infty$. Moreover this sequence satisfies the Armijo Condition \cite[Eqn. 3.4]{Nocedal06}, a popular sufficient decrease condition.
$\blacksquare$\end{example}

There also exist ``second-order'' algorithms that converge only to points at which the Hessian is positive semi-definite \cite{Conn00}. Such algorithms make use of {\em feasible directions of negative curvature}$-$vectors $\vec{d}\in \set{C}$ satisfying $\vec{d}^\top\vec{Hd} < 0$$-$to promote converge to SONC points. Any such algorithm must (periodically) compute a direction of negative curvature and thus, by definition, contains a check of the SOSC: should no direction be found, the SOSC must hold. While this is certainly sufficient when the Hessian is known explicitly, many practical applications do not have such knowledge. Implementations of second-order algorithms that rely on quasi-Newton approximations to the Hessian only determine whether there is a direction of negative curvature for the approximation, rather than the true Hessian, and thus cannot by themselves verify the SOSC. Many other large-scale Hessian-free codes apply Conjugate-Gradient (CG) type iterations to solve constrained quadratic subproblems \cite{Gould01}. Section \ref{SUBSEC:PCG} below, however, demonstrates that convergence of CG methods alone is insufficient to verify the SOSC and thus it is conceivable that CG methods could miss indefiniteness in $\vec{H}$ (over $\set{C}$) in some exceptional circumstances. A post-convergence verification of the SOSC at computed points appears to be required for proper application of Hessian-free methods for large-scale optimization. 


\subsection{Equilibrium Problems}
\label{SUBSEC:EQL}

Verifying the SOSC is also currently vital for solving nonlinear, non-convex equality constrained (generalized Nash) {\em equilibrium} problems; that is, collections of $K \in \{2,3,\dotsc\}$ {\em coupled} optimization problems: 
\begin{equation}
	\label{EQN:Eql}
	\begin{aligned}
		\mnm 	&\quad f_k(\vec{x}_1,\dotsc,\vec{x}_K) \\
		\wrt 		&\quad \vec{x}_k \in \R^{N_k} \\
		\sto 		&\quad \vec{c}_k(\vec{x}_1,\dotsc,\vec{x}_K) = \vec{0}
	\end{aligned}
\end{equation}
where $f_k: \R^{N} \to \R$, $\vec{c}_k : \R^{N} \to \R^{M_k}$, $N_k,M_k \in \N, \; M_k < N_k$, and $N = N_1 + \dotsb + N_K$ for all $k \in \{1,\dotsc,K\}$. A ({\em local}) {\em equilibrium} is a point $(\vec{x}_1^*,\dotsc,\vec{x}_K^*) \in \R^N$ such that $\vec{x}_k^* \in \R^{N_k}$ (locally) solves Problem (\ref{EQN:Eql}) for all $k$. Originating in game theory, such equilibrium models have been used by Economists and Operations Researchers to study electric power \cite{Hobbs07, Hu07, Hu07b, Ehrenmann09} and other energy sectors \cite{Gabriel01, Gabriel05}, new and used vehicles \cite{Berry95, Goldberg95, Goldberg98, Sudhir01, Petrin02, Kleit04, Berry04a, Bento09, Jacobsen10}, entertainment goods \cite{Goolsbee04}, and food services \cite{Besanko98, Nevo00a, Nevo01, Smith04, Thomadsen05}. Recent ``market-systems'' research in engineering design is also applying the equilibrium framework \cite{Michalek04, Skerlos05, Shiau09a, Shiau09b, Shiau09c, Frischknecht10}. See \cite{Facchinei09} for a review of similar applications and economic history. 

First and second order optimality conditions can be employed in the computation of equilibria. In particular, a FONC for a (local) equilibrium follows from combining the FONC for each underlying optimization problem, resulting in a single nonlinear system (e.g., \cite{Morrow10}) or, more generally, Nonlinear Complementarity Problem (NCP) (e.g., \cite{Ferris97, Hu07}); see \cite{Facchinei09}. ``Simultaneously stationary'' points that solve combined FONC can often be computed using standard methods for nonlinear systems or NCPs such as trust-region Newton methods \cite{Dennis96, Conn00}, non-smooth Newton methods \cite{Ralph94, Dirkse95, Munson00}, or semi-smooth methods \cite{Munson00}; see, e.g. \cite{Benson06, Hu07, Morrow10}. However the sufficient decrease conditions that enforce convergence to minimizers or saddle points in optimization problems now apply to a residual norm instead of an objective function, and thus cannot preferentially select equilibria over non-equilibrium stationary points. 

\begin{example}
	\label{EX:EX1}
	Morrow \& Skerlos \cite[Example 10]{Morrow10} compute equilibrium prices for a two-firm market with heterogeneous consumers. Both firms set the price for a single ``branded'' product whose only non-price attribute is ``brand''. There are three types of consumers, two of which are brand-loyal and the other is brand-indifferent. Demand within each type is modeled using a Logit model \cite[Chapter 3]{Train03}. The resulting optimal pricing problem is an unconstrained nonlinear optimization problem with a non-convex, multi-modal objective. The combined FONC for equilibrium has nine solutions, with only four of these nine first-order points locally maximizing both firms' profits. Newton's method applied to the combined FONC cannot distinguish between any of these nine points; any of the five spurious ``solutions'' could be mistaken for equilibria if the SOSC were not verified. Again, no general, globally-convergent method for computing only equilibria in this type of non-convex game is currently known. 
$\blacksquare$\end{example}

Some studies compute equilibria using ``sequential optimization'', ``tattonement'', or ``Jacobi/Gauss-Seidel'' methods \cite{Michalek04, Hu07}; see the discussion of ``practitioners methods'' in \cite{Facchinei09} for algorithmic details. In general, sequential optimization should enable some ``filtering'' of simultaneously stationary points, as the use of optimization algorithms would avoid convergence to simultaneously stationary points that are minimizers of some objectives and maximizers of others. However, there are no existing results guaranteeing {\em any} sort of convergence behavior from sequential optimization methods in either convex or non-convex equilibrium problems \cite{Facchinei09}. Because there do exist algorithms proven to converge to solutions of the combined FONC, computing equilibria through solving the combined FONC is currently theoretically preferable to sequential optimization methods. Furthermore, sequential optimization is likely to be efficient only when the optimization problems are weakly coupled; much effort could be wasted when the optimizers for the coupled problems strongly depend on one another. 

While strong methods for computing equilibrium points in games with {\em convex} objectives and feasible sets exist and are preferable to solving the FONC alone, the alternatives to solving the combined first-order conditions when the players' objectives and feasible sets are non-convex are currently limited. Until methods guaranteed to compute equilibria in non-convex games are developed, general equilibrium programming must be undertaken with checks of the SOSC. 


\subsection{Benefits of Hessian-Free Algorithms}
\label{SUBSEC:HF}

A ``Hessian-free'' algorithm for checking the SOSC will require only matrix-vector products with $\vec{H}$, rather than requiring $\vec{H}$ explicitly \cite[pg. 170]{Nocedal06}. Hessian-free algorithms for checking the SOSC have at least two major advantages over algorithms that require the Hessians explicitly. 

First, multiplying by $\vec{H}$ can be more efficient than working directly with the elements of $\vec{H}$, even when $\vec{H}$ is known \cite[pg. 244]{Trefethen97}. For example: if $\vec{H} \in \R^{N \times N}$ is diagonal, then $\vec{Hx}$ can be computed in $N$ \texttt{flops} rather than the $2N^2$ \texttt{flops} required for arbitrary $\vec{H} \in \R^{N \times N}$. For very large and sparse $\vec{H}$ the efficiency gained by algorithms that require only multiplications by $\vec{H}$ can reduce computational burden by an {\em order of magnitude}, without adding the significant overhead required by sparse factorization methods to track entries and maintain sparsity. The benefits of this property is well-known and lauded in numerical linear algebra. 

Second, the second-order derivatives of the objective and constraints required to explicitly compute $\vec{H}$ can be challenging to derive, difficult to program, and computationally intensive to implement for complex optimization and equilibrium problems \cite{Nash98, Judd98, Nocedal06}. Matrix-vector products $\vec{Hs}$, however, can be obtained with the gradient of the Lagrangian alone using finite-difference approximations \cite{Brown90, Pernice98, Nocedal06}:
\begin{equation*}
	\vec{Hs} \approx 
				\frac{1}{ \sigma \norm{\vec{s}} } 
					\Big( \nabla \set{L} (\vec{x}+\sigma\vec{s},\bsym{\lambda})
								- \nabla \set{L} (\vec{x},\bsym{\lambda}) \Big)
\end{equation*}
for small $\sigma$; see \cite{Dennis96,Brown90} to choose effective scales $\sigma$. In fact, this relationship underlies the effectiveness of some Newton-Krylov solvers for very large and complex nonlinear systems \cite{Brown90, Pernice98}. 

Of course, an approximation to $\vec{H}$ itself could be constructed using finite differences with $N$ gradient evaluations. Indeed, this complete approximation is required to use the BHT and Inertia tests when the second derivatives are not explicitly provided; quasi-Newton approximations such as the BFGS approximation cannot be used. The new Hessian-free algorithms presented below, however, require at most $L \leq N$ evaluations of the gradient of the Lagrangian. For highly constrained problems ($L \ll N$) this represents a significant decrease in function evaluations and thus overall computational burden. 


\section{THREE TESTS OF THE SOSC}
\label{SEC:Tests}

Verifying the SOSC for unconstrained problems requires verifying the positive-definiteness of the Hessian matrix $\vec{H}$. Computationally, this is best accomplished by attempting to take a Cholesky factorization of $\vec{H}$, a stable, efficient, and symmetry-exploiting algorithm \cite{Golub96, Trefethen97, Higham02}. At least three equivalent tests exist for evaluating the ``constrained positive-definiteness'' $\vec{H} \succ_{\set{C}} \vec{0}$ of the Hessian of the Lagrangian in Problem (\ref{EQN:Opt}). 

The most direct test of the SOSC $\vec{H} \succ_{\set{C}} \vec{0}$ is to verify the positive-definiteness of a reduced $L \times L$ Hessian matrix, where $L = N - M$:
\begin{test}[({\em Reduced Matrix Test} \cite{Luenberger09})]
	\label{TEST:CPD}
	$\vec{H} \succ_{\set{C}} \vec{0}$ if, and only if, $\vec{W}^\top\vec{H}\vec{W} \in \R^{L \times L}$ is positive-definite for any matrix $\vec{W} \in \R^{N \times L}$ whose columns form a basis of $\set{C}$.
\end{test}
The three algorithms described in Sections \ref{SUBSEC:CHOL}, \ref{SUBSEC:DIA}, and \ref{SUBSEC:PCG} below are implementations of this test. 

The remaining two tests involve two $(N+M)\times(N+M)$ matrices: 
\begin{equation*}
	\vec{B} = \begin{bmatrix} \vec{0} & \vec{A} \\ \vec{A}^\top & \vec{H} \end{bmatrix}
	\quad\quad\text{and}\quad\quad
	\vec{K} =  \begin{bmatrix} \vec{H} & \vec{A}^\top \\ \vec{A} & \vec{0} \end{bmatrix}
\end{equation*}
where in both cases $\vec{0} \in \R^{M \times M}$. $\vec{B}$ is the well-known ``Bordered Hessian'' \cite{McFadden78, Luenberger09, Papalambros00}. $\vec{K}$ appears in the FONC for equality-constrained quadratic programs and is thus often referred to as a ``KKT matrix'' \cite[Chapter 16]{Nocedal06}. More generally, $\vec{K}$ is an example of a {\em saddle-point} matrix or {\em equilibrium system}; for general information see the extensive review in \cite{Benzi05}. 

The ``Bordered Hessian Test'' (BHT) is a classical test of $\vec{H} \succ_{\set{C}} \vec{0}$ particularly popular in economics that uses $\vec{B}$:
\begin{test}[({\em Bordered Hessian Test} \cite{Gillespie51,McFadden78,Luenberger09})]
	\label{TEST:BHT}
	$\vec{H} \succ_{\set{C}} \vec{0}$ if, and only if, the last $L$ leading principle minors of $\vec{B}$ all have sign $(-1)^M$; specifically, $\mathrm{sign}(\det\big(\vec{B}_i\big)) = (-1)^M$ for all $i = 1,\dotsc,L$, where $\vec{B}_i$ is the submatrix of $\vec{B}$ formed by taking the first $2M+i$ rows and columns. 
\end{test}
Note that the BHT requires $L$ determinant calculations, and thus $L$ LU-factorizations. Section \ref{SUBSEC:BHT} below outlines an efficient procedure for verifying the BHT using updated LU factorization. 

A symmetry-exploiting, single-factorization test based on the KKT matrix $\vec{K}$ can also be derived. The {\em inertia} of a matrix is a triple containing the number of positive, negative, and zero eigenvalues \cite{Golub96, Nocedal06}. Gould \cite{Gould85} proves the following equation:
\begin{lemma}
	\label{LEM:Gould}
	(\cite{Gould85}) $\mathrm{inertia}(\vec{K}) = \mathrm{inertia}(\vec{W}^\top\vec{H}\vec{W}) + (M,M,0)$ for any basis $\vec{W}$ of $\set{C}$. 
\end{lemma}
Because $\vec{H}$ is positive-definite over $\set{C}$ if, and only if, $\mathrm{inertia}(\vec{W}^\top\vec{H}\vec{W}) = (L,0,0)$ $-$ i.e., all of $(\vec{W}^\top\vec{H}\vec{W})$'s eigenvalues are positive $-$ Lemma \ref{LEM:Gould} establishes the following test:
\begin{test}[({\em The Inertia Test} \cite{Gould85})]
	\label{TEST:LDL}
	$\vec{H} \succ_{\set{C}} \vec{0}$ if, and only if, $\mathrm{inertia}(\vec{K}) = \mathrm{inertia}(\vec{B}) = (N,M,0)$.
\end{test}
Section \ref{SUBSEC:LDL} describes how this test can be implemented using an inertia-revealing symmetric-indefinite factorization to compute the inertia of $\vec{K}$ \cite{Bunch77, Gould85, Ashcraft98, Duff04}. 


\section{FIVE ALGORITHMS}
\label{SEC:Algorithms}

This section derives several algorithms for verifying the SOSC; see Table \ref{TAB:Methods}. The focus is on three new Hessian-free algorithms for Test \ref{TEST:CPD}, (Sections \ref{SUBSEC:CHOL}, \ref{SUBSEC:DIA}, and \ref{SUBSEC:PCG}) rather than the existing factorization-based BHT and Inertia tests (Section \ref{SUBSEC:BHT} and \ref{SUBSEC:LDL}). Appendix \ref{APP:Basis} discusses ways to compute a basis $\vec{W}$ of $\set{C}$, needed for the algorithms described in Sections \ref{SUBSEC:CHOL} and \ref{SUBSEC:DIA}; Algorithm \ref{ALG:PCG} provides a Hessian-free method that does not require computing a basis of $\set{C}$. 

	\begin{table}
		\begin{small}
		\caption{Several methods for verifying the SOSC elaborated on in Section \ref{SEC:Algorithms}. ``Hessian-Free'' refers to the ability of an algorithm to operate with Hessian-vector products, rather than the actual elements of the Hessian. ``Basis of $\set{C}$?'' refers to the requirement that an algorithm start by computing a basis of $\set{C}$. ``Dir. of Neg. Curvature'' refers to the ability of an algorithm to reveal or compute a feasible direction of negative curvature when the SOSC fails.}
		\label{TAB:Methods}
		\begin{center}
		\begin{tabular}{ccccc}
				& & {\bfseries Hessian-} 
				& {\bfseries Basis}
				& {\bfseries Dir. of Neg. } \\
			{\bfseries Test}
				& {\bfseries Algorithm} 
				& {\bfseries Free?} 
				& {\bfseries of $\set{C}$?}
				& {\bfseries Curvature} \\ \hline
			\ref{TEST:CPD}
				& Alg. \ref{ALG:CHL} (Sec. \ref{SUBSEC:CHOL})
				& Yes 
				& Yes
				& Yes \\ \hline
			\ref{TEST:CPD}
				& Alg. \ref{ALG:DIA} (Sec. \ref{SUBSEC:DIA})
				& Yes 
				& Yes
				& Yes \\ \hline
			\ref{TEST:CPD}
				& Alg. \ref{ALG:PCG} (Sec. \ref{SUBSEC:PCG})
				& Yes 
				& No
				&  Yes \\ \hline
			 \ref{TEST:BHT}
				& Sec. \ref{SUBSEC:BHT}
				& No
				& No
				& No \\ \hline
			\ref{TEST:LDL}
				& Sec. \ref{SUBSEC:LDL}
				& No 
				& No
				& No \\ \hline
		\end{tabular}
		\end{center}
		\end{small}
	\end{table}


\subsection{An Implicit, Projected Cholesky Factorization for Test (\ref{TEST:CPD})}
\label{SUBSEC:CHOL}

The most direct way to verify $\vec{H} \succ_{\set{C}} \vec{0}$ given a basis $\vec{W}$ of $\set{C}$ is to attempt to take a Cholesky factorization of $\vec{W}^\top\vec{HW}$. In the $n\ith$ step the Cholesky factorization derived in \cite[Lecture 23]{Trefethen97}, 
\begin{align*}
	\vec{W}^\top\vec{HW} 
		&= \vec{L}_1 \dotsb \vec{L}_{n-1}
		\begin{bmatrix}
			\vec{I}_{n-1} & \vec{0} \\
			\vec{0} & \vec{W}_{n:L}^\top \vec{H}_n \vec{W}_{n:L} 
		\end{bmatrix}
		\vec{L}_{n-1}^\top \dotsb \vec{L}_1^\top
\end{align*}
where
\begin{align*}
	\vec{L}_n = \begin{bmatrix} 
				\; \vec{I}_{n-1} \; & \; \vec{0} \; & \; \vec{0} \; \\
				\; \vec{0} \; & \; \sqrt{\alpha_n} \; & \; \vec{0} \; \\
				\; \vec{0} \; & \; \vec{W}_{n+1:L}^\top\vec{H}_n\vec{w}_n / \sqrt{\alpha_n} \; & \; \vec{I}_{L-n-1} \; \\
			\end{bmatrix}, 
\end{align*}
$\vec{H}_1 = \vec{H}$,  
\begin{equation}
	\label{EQN:HnFormula}
	\vec{H}_n = \vec{H}_{n-1} - \vec{H}_{n-1}\vec{w}_{n-1}\vec{w}_{n-1}^\top\vec{H}_{n-1} / \alpha_n, 
\end{equation}
and $\alpha_n = \vec{w}_n^\top\vec{H}_n\vec{w}_n$ for all $n \in \{1,\dotsc,L\}$. $\vec{H} \succ_{\set{C}} \vec{0}$ if, and only if, $\alpha_1,\dotsc,\alpha_L > 0$. 

Observing that the Cholesky factors $\vec{L}_k$ do not need to be explicitly computed to compute the $\alpha$ numbers determining whether the Cholesky process succeeds or fails, the Cholesky process can be reduce to the following:
\begin{lemma}
	For any $n \in \{1,\dotsc,L\}$ such that $\alpha_m = \vec{w}_m^\top\vec{v}_m > 0$ for all $m < n$, $\alpha_n = \vec{w}_n^\top\vec{v}_n$ where
	\begin{align}
		\label{EQN:CHLSum}
		\vec{v}_n
		= \vec{H}_n\vec{w}_n
		=\vec{Hw}_n - \sum_{m=1}^{n-1} \left( \frac{\vec{v}_m^\top\vec{w}_n }{ \alpha_m } \right) \vec{v}_m. 
	\end{align}
\end{lemma}
\proof{}
	The formula for $\alpha_n$ is its definition accepting $\vec{v}_n = \vec{H}_n\vec{w}_n$; the second formula for $\vec{v}_n$ follows from Eqn. (\ref{EQN:HnFormula}) by induction.
\endproof

\begin{figure*}
	\begin{small}
	\begin{minipage}{2.25in}
		\begin{algorithm}[\begin{tiny}(``Classical'' Implicit Cholesky)\end{tiny}]
			\label{ALG:ClassicalCHL}
			\begin{equation*}
				\begin{aligned}
					1&\quad \vec{V} \leftarrow \vec{HW} \\
					2&\quad \texttt{for } n = 1,\dotsc,L, \\
					3&\quad \quad\quad \texttt{for } m = 1,\dotsc,n-1, \\
					4&\quad \quad\quad \quad\quad \beta \leftarrow \vec{v}_m^\top\vec{w}_n / \alpha_m \\
					5&\quad \quad\quad \quad\quad \vec{v}_n \leftarrow \vec{v}_n - \beta \vec{v}_m \\
					6&\quad \quad\quad \alpha_n \leftarrow \vec{w}_n^\top\vec{v}_n \\
					7&\quad \quad\quad \texttt{if }\alpha_n \; \leq \; 0, \; \texttt{fail}, \\
				\end{aligned}
			\end{equation*}
		\end{algorithm}
	\end{minipage}
	\hfill
	\begin{minipage}{2.25in}
		\begin{algorithm}[\begin{tiny}(``Modified'' Implicit Cholesky)\end{tiny}]
			\label{ALG:CHL}
			\begin{equation*}
				\begin{aligned}
					1&\quad \vec{V} \leftarrow \vec{HW} \\
					2&\quad \texttt{for } n = 1,\dotsc,L, \\
					3&\quad \quad\quad \alpha \leftarrow \vec{w}_n^\top\vec{v}_n \\
					4&\quad \quad\quad \texttt{if }\alpha \; \leq \; 0, \; \texttt{fail}, \\
					5&\quad \quad\quad \texttt{for } m = n+1,\dotsc,L, \\
					6&\quad \quad\quad \quad\quad \beta \leftarrow \vec{v}_n^\top\vec{w}_m / \alpha \\
					7&\quad \quad\quad \quad\quad \vec{v}_m \leftarrow \vec{v}_m - \beta \vec{v}_n \\
				\end{aligned}
			\end{equation*}
		\end{algorithm}
	\end{minipage}
	\end{small}
	\caption{``Classical'' (Left) and ``Modified'' (Right) Gram-Schmidt style Implicit Cholesky Algorithms.}
\end{figure*}

Two ``implicit'' Cholesky algorithms that implement Eqn. (\ref{EQN:CHLSum}) are given in Algs. \ref{ALG:ClassicalCHL} and \ref{ALG:CHL}. While $\vec{W}^\top\vec{HW}$ is not explicitly formed, Algs. \ref{ALG:ClassicalCHL} and \ref{ALG:CHL} implicitly form the upper (or lower) triangle of $\vec{W}^\top\vec{HW}$. Note that Algs. \ref{ALG:ClassicalCHL} and \ref{ALG:CHL} only need to compute the products $\vec{Hw}_n$, rather than work with the elements of $\vec{H}$ explicitly, and is thus Hessian-free. Finally, feasible directions of negative curvature can be computed easily should either algorithm reject the SOSC:
\begin{lemma}
	Suppose Algorithm \ref{ALG:ClassicalCHL} or \ref{ALG:CHL} fails in the $n\ith$ step, for some $n \in \{1,\dotsc,L\}$, with $\alpha_n < 0$. Then $\vec{d} = s_1 \vec{w}_1 + \dotsb + s_n \vec{w}_n$ is a feasible direction of negative curvature, where $s_n = 1$ and 
	\begin{equation}
		\label{EQN:CHLs}
		s_m = \frac{ s_{m+1} (\vec{v}_m^\top\vec{w}_{m+1} )
				+ \dotsb + s_n (\vec{v}_m^\top\vec{w}_n) }{\alpha_m} 
	\end{equation}
	for all $m \in \{1,\dotsc,n-1\}$. 
\end{lemma}
Implementing this formula requires storing the $\alpha$ values computed in Alg. \ref{ALG:CHL} and would benefit from storing the triangular array of inner products $\vec{v}_m^\top\vec{w}_{m+k}$ that are computed as part of the Implicit Cholesky process, rather than re-computing them to find a direction of negative curvature. 
\proof{} The idea is straightforward: if $\vec{d} = \vec{Ws}$ where $\vec{s}$ solves $\vec{L}_{n-1}^\top \dotsb \vec{L}_{1}^\top\vec{s} = \vec{e}_n$, then $\vec{d} \in \mathrm{range}(\vec{W}) = \set{C}$ and $\vec{d}^\top\vec{Hd} = \alpha_n < 0$. Eqn. (\ref{EQN:CHLs}) follows from applying back substitution to solve for this $\vec{s}$.
\endproof


\subsection{A Diagonalization Method for Test (\ref{TEST:CPD})}
\label{SUBSEC:DIA}

The following facts furnish a different version of Test (\ref{TEST:CPD}): 
\begin{lemma}
	\label{LEM:TransformationLemma} (i) Let the columns of $\vec{W} \in \R^{N \times L}$ be any basis for $\set{C}$, let $\vec{S} \in \R^{L \times L}$ be any nonsingular matrix, and set $\vec{V} = \vec{WS}$. $\vec{W}^\top\vec{H}\vec{W}$ is positive-definite if, and only if, $\vec{V}^\top\vec{H}\vec{V}$ is positive-definite. (ii) Let $\vec{V}$ be any basis of $\set{C}$ such that $\vec{V}^\top\vec{H}\vec{V}$ is diagonal. $\vec{H} \succ_{\set{C}} \vec{0}$ if, and only if, all of $\vec{V}^\top\vec{H}\vec{V}$'s diagonal entries are positive. 
\end{lemma}
\proof{}
	Claim (i) follows from Sylvester's Law of Inertia and claim (ii) is trivial.
\endproof

Thus, given any basis $\vec{W}$ of $\set{C}$, the definiteness of $\vec{H}$ over $\set{C}$ can be revealed by finding a nonsingular $\vec{S}$ such that $\vec{V} = \vec{WS}$ and $\vec{V}^\top\vec{HV}$ is diagonal. An ``oblique'' Gram-Schmidt process makes this possible: 
\begin{lemma}
	Set $\vec{v}_1 = \vec{w}_1$ and recursively define
	\begin{equation}
		\label{EQN:EGS}
		\vec{v}_n = \vec{w}_n - \sum_{k=1}^{n-1} \left( \frac{\vec{v}_k^\top\vec{Hw}_n}{\alpha_k} \right) \vec{v}_k
	\end{equation}
	for any $n \in \{ 2,\dotsc,L\}$, so long as $\alpha_k = \vec{v}_k^\top\vec{Hv}_k \neq 0$ for all $k \in \{1,\dotsc,n-1\}$. For $\vec{v}_1,\dotsc,\vec{v}_n$ so defined, $\vec{v}_m^\top\vec{Hv}_k = 0$ for all $m,k \in \{1,\dotsc,n\}$, $m \neq k$. 
\end{lemma}
\proof{}
	The proof is a straightforward induction.
\endproof

\begin{figure*}
	\begin{small}
	\begin{minipage}{2.25in}
		\begin{algorithm}[\begin{tiny}(``Classical'' Diagonalization)\end{tiny}]
			\label{ALG:ClassicalDIA}
			\begin{equation*}
				\begin{aligned}
					1&\quad \vec{V} \leftarrow \vec{W} \\
					2&\quad \texttt{for } n = 1,\dotsc,L, \\
					3&\quad \quad\quad \texttt{for } m = 1,\dotsc,n-1, \\
					4&\quad \quad\quad \quad\quad \beta \leftarrow \vec{z}_m^\top\vec{v}_n / \alpha_m \\
					5&\quad \quad\quad \quad\quad \vec{v}_n \leftarrow \vec{v}_n - \beta \vec{v}_m \\
					6&\quad \quad\quad \vec{z}_n \leftarrow \vec{Hv}_n \\
					7&\quad \quad\quad \alpha_n \leftarrow \vec{v}_n^\top\vec{z}_n \\
					8&\quad \quad\quad \texttt{if } \alpha_n \; \leq \; 0, \; \texttt{fail}, \\
				\end{aligned}
			\end{equation*}
		\end{algorithm}
	\end{minipage}
	\hfill
	\begin{minipage}{2.25in}
		\begin{algorithm}[\begin{tiny}(``Modified'' Diagonalization)\end{tiny}]
			\label{ALG:DIA}
			\begin{equation*}
				\begin{aligned}
					1&\quad \vec{V} \leftarrow \vec{W} \\
					2&\quad \texttt{for } n = 1,\dotsc,L, \\
					3&\quad \quad\quad \vec{z} \leftarrow \vec{Hv}_n \\
					4&\quad \quad\quad \alpha_n \leftarrow \vec{v}_n^\top\vec{z} \\
					5&\quad \quad\quad \texttt{if }\alpha_n \; \leq \; 0, \; \texttt{fail}, \\
					6&\quad \quad\quad \texttt{for } m = n+1,\dotsc,L, \\
					7&\quad \quad\quad \quad\quad \beta \leftarrow \vec{w}_m^\top\vec{z} / \alpha_m \\
					8&\quad \quad\quad \quad\quad \vec{v}_m \leftarrow \vec{v}_m - \beta \vec{v}_n \\
				\end{aligned}
			\end{equation*}
		\end{algorithm}
	\end{minipage}
	\end{small}
	\caption{``Classical'' (Left) and ``Modified'' (Right) Gram-Schmidt style Diagonalization algorithms. Note that the vectors $\vec{z}_1,\dotsc,\vec{z}_L$ can overwrite $\vec{w}_1,\dotsc,\vec{w}_L$.}
\end{figure*}

Again, two algorithms for Test \ref{TEST:CPD} based on Eqn. (\ref{EQN:EGS}) are provided in Algs. \ref{ALG:ClassicalDIA} and \ref{ALG:DIA}. Note that Algs. \ref{ALG:ClassicalDIA} and \ref{ALG:DIA} require {\em exactly} the same linear-algebraic operations as Algs. \ref{ALG:ClassicalCHL} and \ref{ALG:CHL}, but on different quantities. In particular, Algs. \ref{ALG:ClassicalDIA} and \ref{ALG:DIA} are also Hessian-free. As written, Algs. \ref{ALG:ClassicalDIA} and \ref{ALG:DIA} require an additional $N \times L$ matrix of storage for $\vec{z}_1,\dotsc,\vec{z}_L$, although these vectors can be written over $\vec{w}_1,\dotsc,\vec{w}_L$ if the basis of $\set{C}$ is not needed after the SOSC check. Note, however, {\em no} additional computation is required to extract a feasible direction of negative curvature from Algs. \ref{ALG:ClassicalDIA} or \ref{ALG:DIA} when the SOSC fails: 
\begin{lemma} 
	Suppose Algorithm \ref{ALG:ClassicalDIA} or \ref{ALG:DIA} fails in step $n \leq L$ with $\alpha < 0$. Then $\vec{v}_n$ is a feasible direction of negative curvature. 
\end{lemma}
Thus, Alg. \ref{ALG:ClassicalDIA} or \ref{ALG:DIA} is a useful re-organization of Alg. \ref{ALG:ClassicalCHL} or \ref{ALG:CHL} (respectively) if directions of negative curvature are important. 


\subsection{Projected Conjugate Gradients for Test (\ref{TEST:CPD})}
\label{SUBSEC:PCG}

The Conjugate Gradient (CG) algorithm is a widely-used iterative method for solving symmetric positive-definite linear systems \cite{Golub96,Trefethen97}. The {\em Projected} Conjugate Gradient (PCG) algorithm is an equivalent algorithm for solving the FONC for equality-constrained quadratic programming problems \cite{Gould01,Nocedal06}. Verifying the SOSC with PCG is based on the following {\em converse} question:
\begin{quote}
	Can (P)CG {\em verify} that a matrix is (constrained) positive-definite?
\end{quote}

This section describes a Hessian-free approach for checking Test (\ref{TEST:CPD}) based on existing PCG algorithms. For simplicity, the majority of the derivation neglects constraints, and considers how CG can be adapted to verify or reject the positive-definiteness of a matrix $\vec{H} \in \R^{N \times N}$. The extension to the constrained case is a straightforward adaptation of this discussion and existing PCG methods as described in, e.g. \cite[Chapter 16]{Nocedal06}. CG applied to the system $\vec{Hx} = \vec{b}$, $\vec{b} \neq \vec{0}$ is denoted by CG$(\vec{H},\vec{b})$; see \cite[pg. 527]{Golub96} or \cite[pg. 112]{Nocedal06} for a formal algorithm. Because CG started at $\vec{x}_0 \neq \vec{0}$ is equivalent to CG$(\vec{H},\vec{b}^\prime)$ started at $\vec{0}$, where $\vec{b}^\prime = \vec{b}-\vec{Hx}_0$, it is assumed that CG$(\vec{H},\vec{b})$ always starts at $\vec{0}$. 

While CG$(\vec{H},\vec{b})$ converges if $\vec{H}$ is definite \cite{Golub96,Trefethen97,Nocedal06}, {\em convergence of CG$(\vec{H},\vec{b})$ alone is not an indicator of definiteness}:
\begin{example} 
	Let $\vec{H} = \mathrm{diag}(1,-1)$ and $\vec{b} = (1,2)$. CG$(\vec{H},\vec{b})$ converges in two steps, despite being applied to an indefinite matrix. However if $\vec{b} = (1,1)$, then CG$(\vec{H},\vec{b})$ breaks down in the first step. 
$\blacksquare$\end{example}
This example shows that convergence (or divergence) of CG alone cannot be used as a basis for verifying positive-definiteness. 

The foundation for using CG$(\vec{H},\vec{b})$ as a positive-definiteness check is, in fact, Lemma \ref{LEM:TransformationLemma}(ii). In exact arithmetic CG$(\vec{H},\vec{b})$ constructs a set of $N^\prime \leq N$ linearly independent, $\vec{H}$-conjugate vectors $\vec{p}_i$ and evaluates $\vec{p}_i^\top\vec{Hp}_i$ for all $i \in \{1,\dotsc,N^\prime\}$ \cite{Golub96, Trefethen97}. Verifying that $\vec{p}_i^\top\vec{Hp}_i > 0$ for all $i$ verifies positive-definiteness, by Lemma \ref{LEM:TransformationLemma}(ii), assuming that $N^\prime = N$. 

\begin{example} Again let $\vec{H} = \mathrm{diag}(1,-1)$ and $\vec{b} = (1,2)$. CG$(\vec{H},\vec{b})$ encounters $\vec{p}_2^\top\vec{Hp}_2 < 0$ in the second step, and thus discovers indefiniteness in $\vec{H}$. $\blacksquare$\end{example}
Ensuring that $\vec{p}_i^\top\vec{Hp}_i > 0$ for all $i \in \{1,\dotsc,N^\prime\}$ is already a component of some CG codes \cite{Nocedal06}. 

Unfortunately CG$(\vec{H},\vec{b})$ may converge for some $N^\prime < N$ prior to building a basis of $\R^N$ and thus cannot always determine whether $\vec{H}$ is positive-definite. Indeed {\em this ``early convergence'' is the central benefit of CG for solving large linear systems}. 
\begin{example} \label{EX:temp} Consider again $\vec{H} = \mathrm{diag}(1,-1)$, and take $\vec{b} = (1,0)$. CG$(\vec{H},\vec{b})$ converges in a single step ($1 = N^\prime < N = 2$), and cannot identify that $\vec{H}$ is indefinite.$\blacksquare$\end{example}
In fact, in exact arithmetic, CG$(\vec{H},\vec{b})$ {\em always} converges with $N^\prime < N$ when $\vec{b}$ lies in a proper invariant subspace of $\vec{H}$. (A {\em proper invariant subspace} of $\vec{H}$ is a proper subspace $\set{W}$ of $\R^N$ such that $\vec{H}\vec{w} \in \set{W}$ for all $\vec{w} \in \set{W}$.) The same issue presents in finite-precision arithmetic, where CG$(\vec{H},\vec{b})$ ``converges'' with $N^\prime < N$ when the residual $\vec{Hx} - \vec{b}$ or the search directions become ``small''.

When CG$(\vec{H},\vec{b})$ converges with $N^\prime < N$, the search can be {\em continued} consistent with Lemma \ref{LEM:TransformationLemma}(ii) by ``restarting'' with a new right-hand-side $\vec{b}^\prime$ that is $\vec{H}$-conjugate to all previous search directions. The continued search should proceed until (a) positive-definiteness has been rejected, (b) the restarted process itself converges, or (c) the remaining dimensions of the space have been searched without rejecting positive-definiteness. 

\begin{example} To continue the case in Example \ref{EX:temp} above choose $\vec{b}^\prime = (0,\pm1)$. $\vec{b}^\prime \perp \vec{Hb}$ and $(\vec{b}^\prime)^\top\vec{H}\vec{b}^\prime < 0$; thus the ``restarted'' CG$(\vec{H},\vec{b}^\prime)$ would identify indefiniteness. 
$\blacksquare$\end{example}

Figure \ref{FIG:PCG} gives a formal algorithm implementing the continued PCG approach$-$with constraints$-$to verify (or reject) $\vec{H} \succ_{\set{C}} \vec{0}$. Algorithm \ref{ALG:PCG} is Hessian-free and, like Algorithm \ref{ALG:DIA}, does not require any additional computation to extract a feasible direction of negative curvature:
\begin{lemma}
	Suppose Algorithm \ref{ALG:PCG} fails in step $j \leq L$ with $\eta < 0$. Then $\vec{p}_j$ is a feasible direction of negative curvature. 
\end{lemma}

Algorithm \ref{ALG:PCG}, being modeled on existing PCG algorithms \cite{Gould01, Nocedal06}, allows flexibility in choosing an operator $\mathrm{proj}_{\set{C}} : \R^N \to \set{C}$ that orthogonally projects vectors from $\R^N$ onto $\set{C} \subset \R^N$. This projection can always be done using a basis $\vec{W}$ of $\set{C}$ by setting $\mathrm{proj}_{\set{C}}(\vec{r}) = \vec{W}(\vec{W}^\top\vec{W})\inv\vec{W}^\top\vec{r}$, a formula that is particularly simple when $\vec{W}$ is orthonormal. However, projection onto $\set{C}$ can also be accomplished without computing a basis for $\set{C}$ \cite{Gould01}. This is an important distinguishing feature of Algorithm \ref{ALG:PCG} relative to Algorithms \ref{ALG:CHL} and \ref{ALG:DIA}: In principle, PCG provides a Hessian-free way to verify $\vec{H} \succ_{\set{C}} \vec{0}$ if it is impractical to project onto $\set{C}$ using a basis of $\set{C}$ computed from the constraint gradients. 

\begin{figure*}
	\begin{small}
	\begin{algorithm}[\begin{tiny}(Continued Projected Conjugate Gradients)\end{tiny}]
		\label{ALG:PCG}
		\begin{equation*}
			\begin{aligned}
				1&\quad \texttt{choose } \vec{b} \in \set{C} \\
				2&\quad i \leftarrow 0 \\
				3&\quad \texttt{while } i \leq L, \\
				4&\quad \quad\quad \texttt{pcg}(\vec{b},\set{C},i) \\
				5&\quad \quad\quad \texttt{if pcg(\vec{b},\set{C},i) failed, exit} \\
				6&\quad \quad\quad \texttt{else if pcg(\vec{b},\set{C},i) converged},\\
				7&\quad \quad\quad\quad\quad \set{C} \leftarrow \set{C} 
											\; \cap \; \mathrm{span}\{ \vec{q}_k : k = i+1,\dotsc,j \}^\perp \\
				8&\quad \quad\quad\quad\quad \texttt{choose } \vec{b} \in \set{C} \\
				9&\quad \quad\quad\quad\quad i \leftarrow j+1 \\
			\end{aligned}
		\end{equation*}
	\end{algorithm}
	\begin{algorithm}[\begin{tiny}(Projected Conjugate Gradient Step)\end{tiny}]
		\begin{equation*}
			\begin{aligned}
				& \texttt{pcg}(\vec{b},\set{C},i) \\
				1&\quad \quad \vec{r} \leftarrow \vec{b}, \;
							\omega \leftarrow \vec{r}^\top\vec{r}, \;
							\vec{r} \leftarrow \vec{r} / \sqrt{\omega}, \;
							\vec{p} \leftarrow \vec{r} \\
				2&\quad \quad \texttt{for } j = i+1,\dotsc,L \\
				3&\quad \quad \quad\quad \tau \leftarrow \omega, \;
									\vec{q}_j \leftarrow \vec{H}\vec{p}, \;
									\eta \leftarrow \vec{p}^\top\vec{q}_j, \\
				4&\quad \quad \quad\quad \texttt{if }  \eta \; \leq \; 0, \texttt{ fail} \\
				5&\quad \quad \quad\quad \vec{r} \leftarrow \vec{r} - (\tau/\eta) \vec{v}_j, \;
									\vec{s} \leftarrow \mathrm{proj}_{\set{C}}(\vec{r}), \;
									\omega \leftarrow \vec{r}^\top\vec{s}, \\
				13&\quad \quad \quad\quad \texttt{if } \abs{\omega} \; \leq \; \texttt{tol}, \texttt{ converged} \\
				14&\quad \quad \quad\quad \vec{p} \leftarrow \vec{s} + (\omega / \tau) \vec{p}
			\end{aligned}
		\end{equation*}
	\end{algorithm}
	\end{small}
	\caption{Formal algorithm for the continued PCG method for verifying (or rejecting) $\vec{H} \succ_{\set{C}} \vec{0}$.}
	\label{FIG:PCG}
\end{figure*}

Gould et al \cite{Gould01} study two projectors in the traditional PCG algorithm (in the context of finite-precision arithmetic). The first sets $\mathrm{proj}_{\set{C}} \vec{r} = \vec{r} - \vec{A}^\top\vec{v}$ where $\vec{AA}^\top\vec{v} = \vec{A}\vec{r}$ can be solved using a Cholesky factorization of $\vec{AA}^\top$. A QR factorization or SVD of $\vec{A}^\top$ could also be used to solve $\min_{\vec{v} \in \R^M} \norm{ \vec{A}^\top\vec{v} - \vec{r} }_2$, the equivalent least-squares problem defining $\vec{v}$. The second method considered in \cite{Gould01} solves an ``augmented system'' with symmetric-indefinite factorization; see \cite{Gould01, Nocedal06, Forsgren07}. Any necessary factorizations need only be taken once, prior to executing Algorithm \ref{ALG:PCG}. Continuation, however, requires projecting onto the {\em intersection} of $\set{C}$ with the image under $\vec{H}$ of the previous search directions. Computationally, this amounts to appending certain rows to $\vec{A}$ or, equivalently, columns to $\vec{A}^\top$. Using Householder QR factorization of $\vec{A}^\top$ provides an efficient and stable way to applying and updating this projector. 


\subsection{Updated LU Factorizations for the BHT, Test (\ref{TEST:BHT})}
\label{SUBSEC:BHT}

The BHT requires computing the sign of $L$ determinants, each typically computed using LU factorization; see \cite[pg. 97]{Golub96} or \cite[pg. 279]{Higham02}. Explicitly forming the $L$ LU factorizations to compute the leading principal minor determinants for the BHT, however, requires an unacceptable and unnecessary amount of work asymptotically proportional to $L(N+M)^3$. 

An efficient LU updating scheme for computing the associated determinant signs follows from a recursive definition of the leading principle submatrices of $\vec{B}$. The matrices $\vec{B}_i$, $i = 2,\dotsc,L$ (defined in Test \ref{TEST:BHT}) satisfy the recursion
\begin{equation*}
	\vec{B}_i
		= \begin{bmatrix}
			 \vec{B}_{i-1} & \vec{b}_i \\
			\vec{b}_i^\top & \gamma_i
		\end{bmatrix}
	\quad\quad\text{for some}\quad\quad
	\begin{aligned}
		&\vec{b}_i \in \R^{2M+i-1} \\
		&\gamma_i \in \R
	\end{aligned}. 
\end{equation*}
An LU factorization of $\vec{B}_i$ can then be constructed by applying the pivots and row eliminations from earlier LU factorizations to $\vec{b}_i$ and then zeroing out $\vec{b}_i^\top$ with $2M+i-1$ new row eliminations \cite[Section 3.2]{Golub96}. 

For brevity the full LU updating process is not described here; deriving a formal statement of this update is straightforward, if a bit tedious. The first LU factorization can be specially designed to account for the leading $M \times M$ submatrix of zeros in $\vec{B}$ and standard pivoting strategies can be used to enhance stability. Existing codes such as \texttt{LAPACK} \cite{LAPACK}, \texttt{LUSOL} \cite{Gill87,LUSOL} \texttt{PARDISO} \cite{Schenk04}, or \texttt{HSL}'s \texttt{LA15} \cite{LA15} can be exploited to take the first LU factorization as well as, in some cases, implement the factorization updates. While computing determinants can be inaccurate due to numerical over- and under-flow \cite[pg. 279]{Higham02}, computing the {\em sign} of the determinant is exact for the computed LU factorization. In any case, the BHT ultimately requires factoring the $L$ matrices $\vec{B}_i$ and may thus be inappropriate for large scale problems without significant sparsity. Furthermore, the BHT does not provide an obvious route to computing a direction of negative curvature should the SOSC fail. 


\subsection{Block Symmetric-Indefinite (LDL) Factorization for the Inertia Test (\ref{TEST:LDL})}
\label{SUBSEC:LDL}

Implementing the Inertia Test (\ref{TEST:LDL}) is, in principle, rather straightforward. The inertia of $\vec{K}$ can be computed using a stable symmetric-indefinite block ``LDL'' factorization \cite{Bunch77, Ashcraft98, Higham02, Duff04} of the form $\vec{P}\vec{K}\vec{P}^\top = \vec{LDL}^\top$ where $\vec{P}$ is a permutation matrix, $\vec{L}$ is unit lower triangular, and $\vec{D}$ is a block-diagonal matrix with $1 \times 1$ and $2 \times 2$ blocks; see \cite[Chapter 11]{Higham02}. By Sylvester's Law of Inertia, $\mathrm{inertia}(\vec{K}) = \mathrm{inertia}(\vec{D})$ \cite[pg. 403]{Golub96} and computing the inertia of the computed $\vec{D}$ is efficient and exact even in inexact arithmetic; see \cite[Problem 11.2]{Higham02}. Indeed, symmetric-indefinite codes \cite{Duff04,Schenk06} already have options to compute inertia, and the matrix inertia is used in some optimization solvers to promote global convergence \cite{Gill91,Byrd06,Schenk08}. However accuracy of the Inertia test in inexact arithmetic depends on whether the inertia of the $\vec{D}$ computed using finite-precision arithmetic actually equals the inertia of $\vec{K}$. The Inertia test requires factoring a single $(N+M) \times (N+M)$ matrix and may thus be applicable to large scale problems only when there is significant sparsity. There are, however, extremely efficient codes for forming this factorization. Like the BHT, it is not obvious how to compute a direction of negative curvature from $\vec{K}$ when the Inertia Test fails. 


\section{EXAMPLES}
\label{SEC:Examples}

This section presents several computational examples. Sections \ref{SUBSEC:GoodA} and \ref{SUBSEC:IllA} illustrate the relative performance characteristics of the different algorithms discussed above on dense test problems for which $\vec{H}$, $\vec{A}$, and the truth of $\vec{H} \succ_{\set{C}} \vec{0}$ is known explicitly. Section \ref{SUBSEC:OITP} demonstrates the advantages of the Hessian-free properties of Cholesky and Diagonalization using a test problem from the \texttt{COPS} collection \cite{Dolan04}. 

\subsection{Computational Details}
\label{SUBSEC:CompDetails}

Algs. \ref{ALG:CHL}, \ref{ALG:DIA}, \ref{ALG:PCG}, the updated LU BHT, and the Inertia test have been implemented in \texttt{C} making extensive use of \texttt{BLAS} \cite{BLAS} routines and \texttt{LAPACK}'s routines for QR (\texttt{dgeqrf}/\texttt{dormqr}), LU (\texttt{dgetrf}), and symmetric-indefinite block LDL (\texttt{dsytrf}) factorizations \cite{LAPACK}. The basis $\vec{W}$ of $\set{C}$ used in Algorithms \ref{ALG:CHL} and \ref{ALG:DIA} is computed using a QR factorization of $\vec{A}^\top$, primarily for consistency with the updated QR approach to the PCG algorithm. Using the SVD would probably provide more stable (though also more expensive) computations of $\vec{W}$. The updated LU BHT described above was verified using a ``naive'' BHT that computes $\mathrm{sign}(\det(\vec{B}_i))$ via $L$ independent LU factorizations. The time savings from updating LU factorizations ranged from 10-90\% depending on $N$ and $M$. All computations were undertaken on an Apple MacPro tower with dual quad-core 2.26 GHz ``Nehalem'' processors (each with an 8 MB cache) and 32 GB of RAM running Mac OS X 10.6.6 and Apple's implementations of \texttt{BLAS} and \texttt{LAPACK}. 

\subsection{Generating Dense Random Test Problems}
\label{SUBSEC:GenRandProb}

Random test problems were obtained as follows: Let $N \in \{2,3,\dotsc\}$, $M \in \{1,\dotsc,N-1\}$, and $P \in \{0,\dotsc,N\}$. Choose symmetric positive-definite $\bsym{\Lambda}_+ \in \R^{P \times P}$, symmetric negative-definite $\bsym{\Lambda}_- \in \R^{(N-P)\times(N-P)}$, orthogonal $\vec{Q} \in \R^{N \times N}$, upper-triangular $\vec{R} \in \R^{M \times M}$, and set 
\begin{equation}
	\label{EQN:HAFormula}
	\vec{H} = \vec{Q} \begin{bmatrix} \bsym{\Lambda}_+ & \vec{0} \\ \vec{0} & \bsym{\Lambda}_- \end{bmatrix} \vec{Q}^\top
	\quad\quad\text{and}\quad\quad
	\vec{A}^\top = \vec{Q} \begin{bmatrix} \vec{0} \\ \vec{R} \end{bmatrix}
\end{equation}
Whether $\vec{H} \succ_{\set{C}} \vec{0}$ for such $\vec{H}$ and $\vec{A}$ is known analytically: 
\begin{lemma}
	\label{LEM:ExampleLem}
	Let $\vec{H}$ and $\vec{A}$ be defined as in Eqn. (\ref{EQN:HAFormula}) and let $\set{C} = \mathrm{null}(\vec{A})$. (i) $\vec{H} \succ_{\set{C}} \vec{0}$ if, and only if, $L = N - M \leq P$. (ii) For uniformly drawn $M \in \{1,\dotsc,N-1\}$ and $P \in \{0,\dotsc,N\}$, $\Prob(\vec{H} \succ_{\set{C}} \vec{0}) = (N+2)/(2(N+1))$. 
\end{lemma}
The techniques described by Higham \cite[pg, 517-518]{Higham02} are used to draw random orthogonal matrices ($\vec{Q}$) and random symmetric matrices with known eigenvalues ($\bsym{\Lambda}_+, \bsym{\Lambda}_-$). Note that any such construction is likely to be dense, and thus does not take advantage of the Hessian-free characteristics of the Cholesky, Diagonalization, or PCG algorithms. 

\subsection{A Well-Conditioned Example}
\label{SUBSEC:GoodA}

\begin{example}
	\label{EX:GoodA}
	Let $\vec{A}$ be defined by an $\vec{R}$ with off diagonal elements $r_{i,j}$, $j > i$, drawn from a standard normal distribution and diagonal elements $r_{i,i}$ drawn from a normal distribution with mean zero and variance $(M-i)^2$. 
$\blacksquare$\end{example}

Over 20,000 numerical tests were run for Ex. \ref{EX:GoodA}. $90$ distinct values of $N$ were drawn from $\{10,...,5000\}$ and, for each $N$, at least $200$ $M,P$ pairs were drawn. Slightly more than half ($\sim$ 50.3\%) of these trials have $\vec{H} \succ_{\set{C}} \vec{0}$, with $\vec{H} \nsucc_{\set{C}} \vec{0}$ in the remaining trials; see Lemma \ref{LEM:ExampleLem}(ii). Every Hessian matrix $\vec{H}$ drawn has eigenvalues $\lambda$ satisfying $0.1 \leq \abs{\lambda} \leq 100$. The tolerance for convergence in the PCG approach (Alg. \ref{ALG:PCG}) was $10^{-10}$. 

In Ex. \ref{EX:GoodA}, all methods correctly verified or rejected $\vec{H} \succ_{\set{C}} \vec{0}$ in most tests. No method returned a false positive in any trial, and only the Inertia test returned false negatives. Indeed, the performance of the Inertia test degrades when $\vec{H} \succ_{\set{C}} \vec{0}$ as $N$ grows, as shown in Table \ref{TAB:MITFailures}. For $N \approx 5,000$, more than 3\% of the SOSC checks using the Inertia test gives a false negative. No other algorithm returned a single false negative. 

\begin{table*}
	\label{TAB:MITFailures}
	\begin{small}
	\caption{Percent of trials in Example \ref{EX:GoodA} for which $\vec{H} \succ_{\set{C}} \vec{0}$ yet the Inertia Test incorrectly declares $\vec{H} \nsucc_{\set{C}} \vec{0}$; i.e. false negatives (FN).}
	\begin{center}
	\begin{tabular}{lcccccccccccccccc}
		$N$ & 77 & 242 & 478 & 679 & 899 & 908 & 1153 & 1161 \\ \hline
		FN & 1.0\% & 0.5\% & 2.0\% & 1.1\% & 1.0\% & 1.0\% & 1.0\% & 2.1\% \\ \hline
		\\ %
		$N$ & 1530 & 1538 & 1718 & 2306 & 2460 & 3687 & 4637 & 4903 \\ \hline
		FN & 1.0\% & 1.0\% & 2.9\% & 2.0\% & 1.1\% & 2.7\% & 1.9\% & 3.3\% \\ \hline
	\end{tabular}
	\end{center}
	\end{small}
\end{table*}

\begin{figure*}
	\includegraphics[width=1.75in]{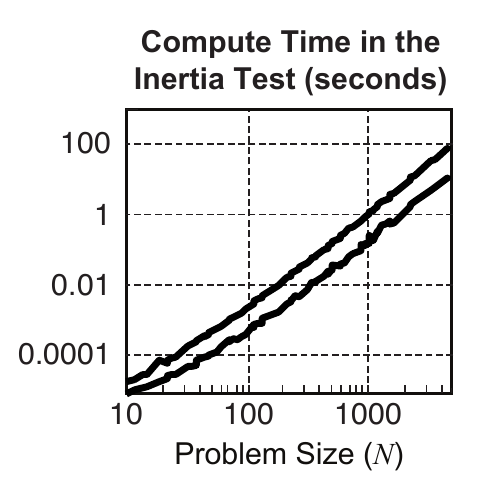}
	\caption{Compute times for the Inertia test applied in Example \ref{EX:GoodA} in trials for which $\vec{H} \succ_{\set{C}} \vec{0}$ is true.}
	\label{FIG:EX1TrueTimes}
\end{figure*}

\begin{figure*}
	\includegraphics[width=5.25in]{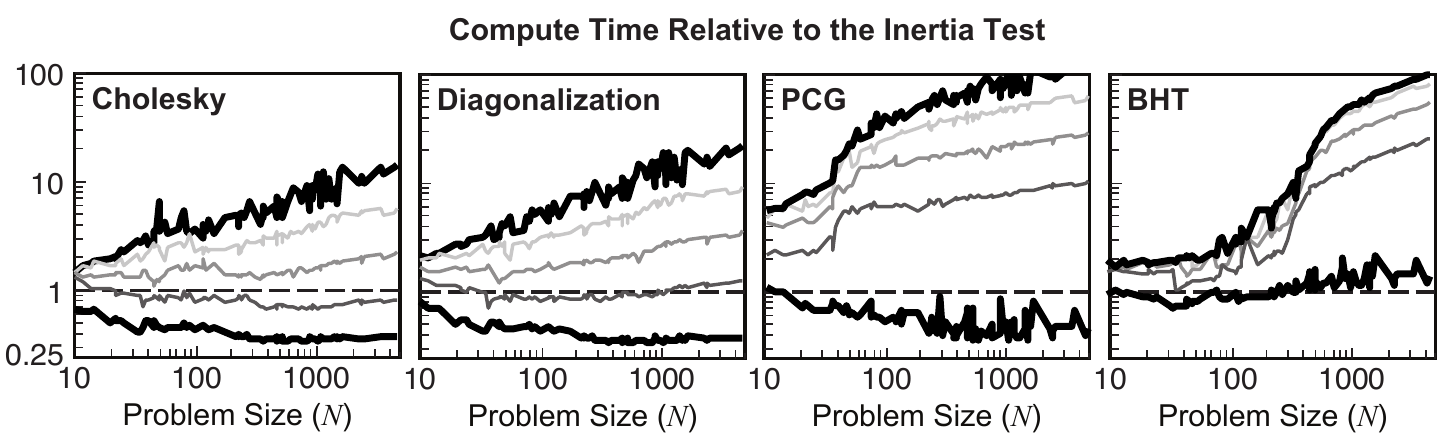}
	\caption{Compute times for Example \ref{EX:GoodA} in trials for which $\vec{H} \succ_{\set{C}} \vec{0}$ is true, relative to the Inertia test. Solid black lines give the maximum and minimum ratios or times for a given $N$ over the values of $M$ tested. Thin grey lines in the four relative plots represent the maximum ratio of times for subsets of trials such that $M \geq 0.25N$ (lightest), $M \geq 0.5N$ (midtone), and $M \geq 0.75N$ (darkest).}
	\label{FIG:EX1RelTimes}
\end{figure*}

Figs. \ref{FIG:EX1TrueTimes} and \ref{FIG:EX1RelTimes} compare the time required by each method to computationally verify $\vec{H} \succ_{\set{C}} \vec{0}$ for Ex. \ref{EX:GoodA}. The Inertia test, generally the fastest method, is used as a benchmark for comparison; see Fig. \ref{FIG:EX1TrueTimes}. Cholesky (Alg. \ref{ALG:CHL}) and Diagonalization (Alg. \ref{ALG:DIA}) are as fast or faster than the Inertia test for highly constrained problems with $M \geq 0.75N$. However if $M$ is small relative to $N$, Algs. \ref{ALG:CHL} and \ref{ALG:DIA} can take 10-20 times longer than the Inertia test. Alg. \ref{ALG:DIA} is slightly slower than Alg. \ref{ALG:CHL} in part because the Hessian-vector products $\vec{Hw}$ are done one-by-one in Alg. \ref{ALG:DIA}, instead of ``all-at-once'' in Alg. \ref{ALG:CHL}. When $\vec{H}$ is known the level-3 \texttt{BLAS} routines used to form $\vec{HW}$ optimize efficient cache memory usage. Generally speaking these results are more encouraging than they appear: \texttt{dsytrf}, the code for the Inertia test, is highly optimized while the implementations of Algs. \ref{ALG:CHL} and \ref{ALG:DIA} have not yet been optimized. ``Blocked'' variants of Algs. \ref{ALG:CHL} and \ref{ALG:DIA} that fully exploit memory traffic efficiencies in the level-3 \texttt{BLAS} are conceptually easy to derive, and will be even more competitive with the Inertia test from the perspective of compute time. 

The PCG approach (Alg. \ref{ALG:PCG}) and the BHT are as fast or faster than the Inertia test only for $M \approx N$, and can take more than 100 times longer than the Inertia test on problems with $N \geq 1000$ and small $M$. Though it is not shown in the plots, continuation is an important part of the PCG approach; see Table \ref{TAB:Continuation}. Alg. \ref{ALG:PCG} converged at least once in over 68\% of the trials for which $\vec{H} \succ_{\set{C}} \vec{0}$, and would thus ``fail'' in more than 68\% of our cases if it were not continued. Moreover as $N$ grows, Alg. \ref{ALG:PCG} tends to converge more often. Based on the trials undertaken for Ex. \ref{EX:GoodA}, we would ``expect'' Alg. \ref{ALG:PCG} to converge in $\sim N^{0.85}$ cases. Specifically, Alg. \ref{ALG:PCG} converged no more than $N^{0.9}$ times as $N$ grew in the trials undertaken for Ex. \ref{EX:GoodA}, and converged at least $\sim N^{0.8}$ times in more than 50\% of the trials. While these predictions should not be extrapolated beyond this example, they clearly demonstrate the necessity of continuation in the PCG approach to verifying the SOSC. 

\begin{table}
	\begin{center}
	\begin{small}
	\label{TAB:Continuation}
	\caption{Percent of trials in Example \ref{EX:GoodA} for which $\vec{H} \succ_{\set{C}} \vec{0}$ and PCG was continued at least once.}
	\begin{tabular}{cccccc}
		All $N$ & $N \geq 20$ & $N \geq 50$ & $N \geq 100$ & $N \geq 1,000$ \\ \hline
		68\% & 76\% & 88\% & 94\% & 99\% \\ \hline
	\end{tabular}
	\end{small}
	\end{center}
\end{table}

\subsection{A Poorly Conditioned Example}
\label{SUBSEC:IllA}

\begin{example}
	\label{EX:IllA}
	Let $\vec{A}$ be defined by an $\vec{R}$ with all elements drawn from a standard normal distribution. 
$\blacksquare$\end{example}

Again, over 20,000 numerical tests were run for Ex. \ref{EX:IllA} with the same character as for Ex. \ref{EX:GoodA}, discussed above. 
	
\begin{figure*}
		\begin{center}
		\includegraphics{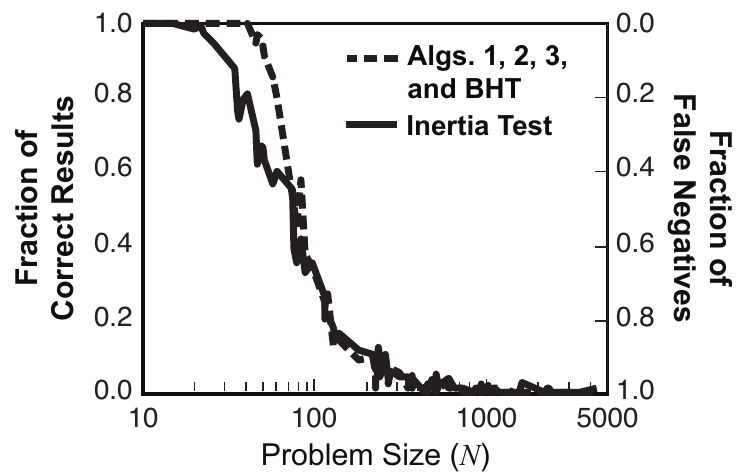}
		\end{center}
		\caption{Fraction of correct test results (left axis) and, symmetrically, false negatives (right axis), over all trials in Example \ref{EX:IllA} for which $\vec{H} \succ_{\set{C}} \vec{0}$. }
		\label{FIG:Correct}
	\end{figure*}

Ex. \ref{EX:IllA} demonstrates that numerical accuracy is far from guaranteed for computational SOSC checks. Fig. \ref{FIG:Correct} illustrates the fraction of correct tests results and thus also of false negatives. The Cholesky, Diagonalization, PCG, and BHT tests appear more stable than the Inertia test. However all methods have a false negative rate close to 70\% for problems with as few as $N \approx 100$ variables. With $N \approx 1,000$ variables, each test is so overcome by roundoff error that virtually no correct results are obtained. As with Ex. \ref{EX:GoodA}, there were no false positives. 

The relatively poorer performance of the Inertia test in both examples deserves some discussion. Weyl's Eigenvalue Pertubation Theorem provides one way to control numerical error in the Inertia test: 
\begin{theorem}
	(Weyl's Theorem \cite[Theorem 4.1]{Eisenstat98}) If $\lambda_1 \geq \dotsb \geq \lambda_n$ and $\mu_1 \geq \dotsb \geq \mu_n$ are the eigenvalues of $\vec{X}$ and $\vec{X}+\vec{E}$, respectively, for any $n \times n$ symmetric matrix $\vec{X}$ and any $n \times n$ matrix $\vec{E}$, then $\abs{ \lambda_i - \mu_i } \leq \norm{\vec{E}}_2$ for all $i$. 
\end{theorem}
\begin{corollary}
	Let $\vec{P}(\vec{K}+\vec{E})\vec{P}^\top = \vec{LDL}^\top$ be the {\em computed} LDL factorization of $\vec{K}$ where errors are accumulated into a perturbation $\vec{E}$ of $\vec{K}$; see \cite[Chapter 11]{Higham02}. If $\abs{\mu}_{\min} = \min\{ \abs{\mu} : \mu \in \lambda(\vec{K}+\vec{E}) \} > \norm{\vec{E}}_2$ then $\mathrm{inertia}(\vec{K}) = \mathrm{inertia}(\vec{D})$. 
\end{corollary}
\proof The Weyl Theorem states that $\lambda_i = \mu_i + \epsilon_i$ for some $\abs{\epsilon_i} \leq \norm{\vec{E}}_2$ for all $i$. If $\abs{\mu_i} > \norm{\vec{E}}_2$ for all $i$, then also $\mathrm{sign}\{\lambda_i\} = \mathrm{sign}\{\mu_i\}$. Hence $\mathrm{inertia}(\vec{K}) = \mathrm{inertia}(\vec{K} + \vec{E}) = \mathrm{inertia}(\vec{D})$.
\endproof
If, however, $\abs{\mu}_{\min} \leq \norm{\vec{E}}_2$ then $\mathrm{inertia}(\vec{K}+\vec{E}) = \mathrm{inertia}(\vec{D})$ may {\em not} be equal to $\mathrm{inertia}(\vec{K})$. That is, if the magnitude of any eigenvalue of $\vec{K}+\vec{E}$ is less than the norm of the accumulated round-off errors, this eigenvalue of the perturbed KKT matrix {\em may} have the wrong sign relative to the true eigenvalue. As $N+M$ grows, the likelihood of both a small $\abs{\mu}_{\min}$ {\em and} large $\norm{\vec{E}}_2$ increases. 

Note that $\vec{K}$ can have small eigenvalues even if $\vec{H}$ does not. Specifically, $\vec{K}$ has small eigenvalues whenever $\vec{A}$ has nearly linearly dependent columns, as the simple example in Appendix \ref{APP:KExample} shows. 
	
	Using the ``exact'' basis $\vec{Q}_{1:L}$ of $\set{C}$ in the Cholesky and Diagonalization tests eliminated the false negatives seen in Fig. \ref{FIG:Correct}. This suggests that the numerical errors in these two tests, as well as perhaps the PCG test, are entirely a consequence of numerical errors in the QR factorization of $\vec{A}^\top$. Two remarks along these lines must be made: First, using $\vec{Q}_{1:L}$ is a device of the artificial numerical examples. In practice, some factorization must be used compute a basis for $\set{C}$ from $\vec{A}$. Second, the effect of error in the computed basis of $\set{C}$ will depend on the spectrum of $\vec{H}$. For example if $\vec{H}$ is positive-definite these representational errors in would be irrelevant to the accuracy of the test. On the other hand, false negatives can be obtained only if the numerical errors in the representation of $\set{C}$ magnify the contribution of $\vec{H}$'s negative eigenvalues to the quadratic forms. Understanding when the construction of a basis of $\set{C}$ is sufficiently accurate relative to $\vec{H}$ will be essential to understanding when SOSC checks using the Cholesky and Diagonalization approaches are themselves numerically accurate. 
		
	Unfortunately there will also be numerical errors associated with the Cholesky and Diagonalization approaches themselves for some Hessians $\vec{H}$, independent of what errors are made in the construction of a representation of $\set{C}$. What specific properties of $\vec{H}$ to monitor and control in this respect are not yet known and will be the subject of future investigations. 

\subsection{A Hessian-Free Example from the \texttt{COPS} Collection}
\label{SUBSEC:OITP}

The ``Thomson Problem'' of finding the minimal energy configuration of $K \in \{ 2,3,\dotsc\}$ points on a sphere is a a large-scale equality-constrained optimization problem from the \texttt{COPS} collection \cite{Dolan04}:
\begin{equation}
	\label{EQN:TP}
	\begin{aligned}
		\mnm &\quad f(\vec{x}_1,\dotsc,\vec{x}_K) 
					= \sum_{k=1}^{K-1} \sum_{m=k+1}^K \frac{1}{\norm{\vec{x}_k - \vec{x}_m}_2} \\
		\wrt &\quad \vec{x}_1,\dotsc,\vec{x}_K \in \R^3 \\
		\sto &\quad \norm{\vec{x}_k}_2^2 = 1 \quad\text{for all}\quad n \in \{1,\dotsc,K\}
	\end{aligned}
\end{equation}
Prob. (\ref{EQN:TP}) depends only on the Euclidean norms of $K$ vectors in $\R^3$, which are invariant over orthogonal transformations of $\R^3$. Thus, Prob. (\ref{EQN:TP}) as stated has infinitely many solutions: specifically, if $(\vec{x}_1,\dotsc,\vec{x}_K)$ is a solution then so is $(\vec{Qx}_1,\dotsc,\vec{Qx}_K)$ for each rotation or reflection $\vec{Q}$ of $\R^3$. An ``orthogonally-invariant'' Thomson problem can be obtained by restricting the locations of the first and second points in a manner inspired by Householder QR factorization: 
\begin{equation}
	\label{EQN:OITP}
	\begin{aligned}
		\mnm &\quad f(\vec{x}_1,\dotsc,\vec{x}_K) 
					= \sum_{k=1}^{K-1} \sum_{m=k+1}^K \frac{1}{\norm{\vec{x}_k - \vec{x}_m}_2} \\
		\wrt &\quad \vec{x}_1,\dotsc,\vec{x}_K \in \R^3 \\
		\sto &\quad \norm{\vec{x}_k}_2^2/2 = 1/2 \quad\text{for all}\quad n \in \{1,\dotsc,K\} \\
			&\quad x_{1,2} = x_{1,3} = x_{2,3} = 0 \\
	\end{aligned}
\end{equation}
Problem (\ref{EQN:OITP}) has $N = 3K$ variables and $M = K + 3$ constraints. 

Prob. (\ref{EQN:OITP}) was solved for specific values of $K$ between 3 and 333 using \texttt{matlab}'s Hessian-free interior-point algorithm. The SOSC was then verified at the computed $(\vec{x}_1,\dotsc,\vec{x}_K)$ and $\bsym{\lambda}$ using our \texttt{C} implementations of Cholesky, Diagonalization, and the Inertia test. The Hessian-free nature of the Cholesky and Diagonalization algorithms was exploited by using directional finite differences to approximate Hessian-vector products $\vec{Hs}$. For the Inertia test the full Hessian $\vec{H}$ was approximated with finite-differences. The PCG and BHT algorithms were not used; the results in Section \ref{SUBSEC:GoodA} suggest they are not competitive. 

\begin{figure*}
		\begin{center}
		\includegraphics[width=2in]{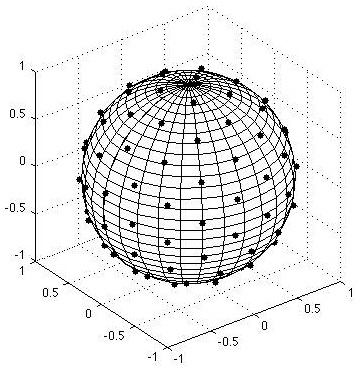}
		$\quad$
		\includegraphics[width=3in]{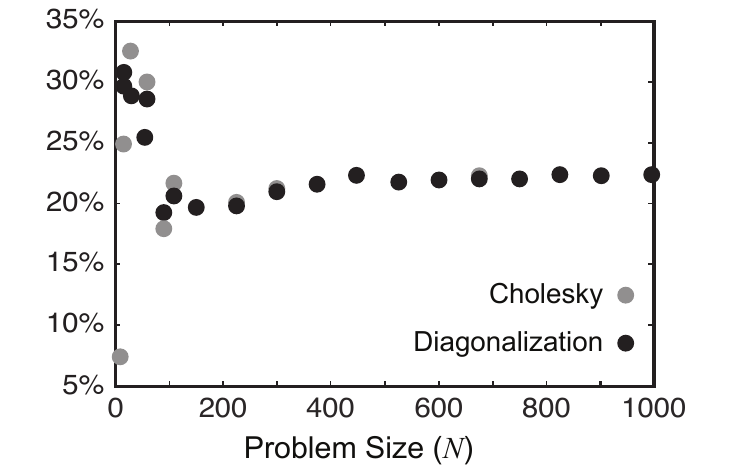}
		\end{center}
		\caption{Left: Computed solution to Prob. (\ref{EQN:OITP}) for $K = 100$. Right: Percent reduction in time to computationally verify $\vec{H} \succ_{\set{C}} \vec{0}$ for various instances of Prob. (\ref{EQN:OITP}). When only the black dot is visible, the gray and black dots coincide. Note that $N = 3K$, where $K$ is the number of points distributed over the surface of the unit sphere. }
		\label{FIG:OITP}
	\end{figure*}

In this case the Hessian-free SOSC checks provided by the Cholesky and Diagonalization algorithms tended to reduced time to verify the SOSC relative to the Inertia test by just over 20\%; see Fig. \ref{FIG:OITP}. Recall that the Inertia test was the fastest method on the dense test problems above with $\sim$ 30\% of the variables constrained. Thus the Hessian-free application of the Cholesky and Diagonalization algorithms results in significant computational savings. 

\section{CONCLUSIONS}

This article has presented three novel Hessian-free algorithms for verifying (or rejecting) the SOSC for constrained continuous optimization. These algorithms also make computation of feasible directions of negative curvature easy when the SOSC fails, a feature not available in classical tests. Numerical trials have demonstrated (1) the inefficiency of the Bordered Hessian Test, (2) the relative speed of the Inertia test, (3) the computational efficiency of the new algorithms, especially when their Hessian-free properties can be exploited, and (4) the potential for significant loss of accuracy due to round-off error using any method, even on small problems. Future work will optimize implementations of the new algorithms and undertake a detailed mathematical analysis of round-off errors to determine computable certificates of test accuracy. 


\bibliographystyle{amsplain}
\bibliography{references}


\appendix       
\section{Computing a Basis for $\set{C}$ given $\vec{A}$}
\label{APP:Basis}

Both Algorithm \ref{ALG:CHL} and \ref{ALG:DIA} require a basis $\vec{W}$ for $\set{C}$, rather than the constraint gradients $\vec{A}$. There are several ways to compute such a basis from $\vec{A}$. Implicitly, each method considers the equation $\vec{AW} = \vec{0}$ and uses some factorization of $\vec{A}$ (or $\vec{A}^\top$) to find a formula for $\vec{W}$. Several examples follow:
\begin{itemize}
	\item SVD of $\vec{A}$: If $\vec{A} = \vec{U}[ \; \bsym{\Sigma} \;\; \vec{0} \; ]\vec{V}^\top$ is a full Singular Value Decomposition (SVD; \cite[Lecture 4]{Trefethen97}) of $\vec{A}$, then $\vec{W} = \vec{V}_{M+1:N}$ is an {\em orthonormal} basis for $\set{C}$. 
	\item QR factorization of $\vec{A}^\top$: Similarly if $\vec{A}^\top = \vec{QR}$ then $\vec{W} = \vec{Q}_{M+1:N}$ is an {\em orthonormal} basis for $\set{C}$. 
	\item QR or LU factorization of $\vec{A}$: If $\vec{A} = \vec{Q}[ \; \vec{R} \;\; \vec{S} \; ]$ for upper-triangular $\vec{R}$, or if $\vec{PA} = \vec{L}[ \; \vec{U} \;\; \vec{S} \; ]$ for upper-triangular $\vec{U}$ and some permutation matrix $\vec{P}$, then 
	\begin{equation*}
		\vec{W} = \begin{bmatrix} \vec{R}\inv\vec{S} \\ \vec{I} \end{bmatrix}
		\quad\quad\text{or}\quad\quad
		\vec{W} = \begin{bmatrix} \vec{U}\inv\vec{S} \\ \vec{I} \end{bmatrix}
	\end{equation*}
	are bases for $\set{C}$. In both cases, $\vec{I} \in \R^{L \times L}$. 
\end{itemize}
Stability considerations typically suggest that either of the first two approaches are preferable to the third. 

\section{$\vec{K}$ has Small Eigenvalues When $\vec{A}$ is Nearly Rank Deficient}
\label{APP:KExample}

Let $\vec{H} \in \R^{N \times N}$ and $\vec{A}^\prime \in \R^{(M-1)\times N}$ be arbitrary (but full-rank). Let $\vec{a}_m \in \R^N$ denote the $m\ith$ {\em row} of $\vec{A}^\prime$, for $m \in \{1,\dotsc,M-1\}$. Choose $\beta_1,\dotsc,\beta_{M-1}$ (not all zero), $\bsym{\epsilon} \in \R^N$ with small 2-norm, and set $\vec{a}_M = \sum_{m=1}^{M-1} \beta_m \vec{a}_m + \bsym{\epsilon}$. The $\beta$'s and $\bsym{\epsilon}$ can be chosen so that $\norm{\vec{a}_M}_2 = 1$, even if $\norm{\bsym{\epsilon}}_2$ is very small. Define
\begin{equation*}
	\vec{K} = \begin{bmatrix} \vec{H} & (\vec{A}^\prime)^\top & \vec{a}_M \\
						\vec{A}^\prime & \vec{0} & \vec{0} \\
						\vec{a}_M^\top & \vec{0} & 0 \end{bmatrix}
	\quad\quad\text{and}\quad\quad
	\vec{F} = \begin{bmatrix} 
				\vec{0} & \vec{0} & - \bsym{\epsilon} \\ 
				\vec{0} & \vec{0} & \vec{0} \\
				\vec{0} & \vec{0} & 0 
			\end{bmatrix}, 
\end{equation*}
noting that $\vec{K} + \vec{F}$ is singular. Specifically, 
\begin{align*}
	\left( \vec{K} + \begin{bmatrix} 
				\vec{0} & \vec{0} & - \bsym{\epsilon} \\ 
				\vec{0} & \vec{0} & \vec{0} \\
				\vec{0} & \vec{0} & 0 
			\end{bmatrix} \right)
		\begin{bmatrix}
			\vec{0} \\ 
			- \bsym{\beta} \\
			0
		\end{bmatrix}
		&= \begin{bmatrix} 
			\vec{H} & (\vec{A}^\prime)^\top & \vec{a}_M - \bsym{\epsilon} \\
			\vec{A}^\prime & \vec{0} & \vec{0} \\
			\vec{a}_M^\top & \vec{0} & 0 
		\end{bmatrix}
		\begin{bmatrix}
			\vec{0} \\ 
			- \bsym{\beta} \\
			0
		\end{bmatrix} \\
		&= \begin{bmatrix}
			- \sum_{m=1}^{M-1} \beta_m \vec{a}_m + \vec{a}_M - \bsym{\epsilon} \\
			\vec{0} \\
			0
		\end{bmatrix}
		= \begin{bmatrix}
			\vec{0} \\
			\vec{0} \\
			0
		\end{bmatrix}. 
\end{align*}
Thus, $\vec{K}+\vec{F}$ has a zero eigenvalue and Weyl's Theorem states that $\vec{K}$ has an eigenvalue $\lambda$ with $\abs{\lambda} \leq \norm{\vec{F}}_2 = \norm{\bsym{\epsilon}}_2$. 


\end{document}

%% file: customcommands.tex
\newcommand{\N}{\ensuremath{\mathbb{N}}}

\newcommand{\R}{\ensuremath{\mathbb{R}}}

\newcommand{\abs}[1]{\ensuremath{\left\lvert #1 \right\rvert}}
\newcommand{\norm}[1]{\ensuremath{\lvert\lvert #1\rvert\rvert}}

\newcommand{\bsym}[1]{\ensuremath{\boldsymbol{#1}}}
\renewcommand{\vec}[1]{\ensuremath{\mathbf{#1}}}
\newcommand{\set}[1]{\ensuremath{\mathcal{#1}}}

\newcommand{\inv}{\ensuremath{^{-1}}}
\newcommand{\Prob}{\ensuremath{\mathbb{P}}}

\newcommand{\mnm}{\ensuremath{\text{minimize}}}
\newcommand{\wrt}{\ensuremath{\mathrm{with \; respect \; to}}}
\newcommand{\sto}{\ensuremath{\mathrm{subject \; to}}}

\newcommand{\ith}{\ensuremath{^\text{th}}}

%% file: SOSC.bbl
\providecommand{\bysame}{\leavevmode\hbox to3em{\hrulefill}\thinspace}
\providecommand{\MR}{\relax\ifhmode\unskip\space\fi MR }
\providecommand{\MRhref}[2]{%
  \href{http://www.ams.org/mathscinet-getitem?mr=#1}{#2}
}
\providecommand{\href}[2]{#2}
\begin{thebibliography}{10}

\bibitem{Ashcraft98}
Cleve Ashcraft, Roger~G. Grimes, and John~G. Lewis, \emph{Accurate symmetric
  indefinite linear equation solvers}, SIAM journal of matrix analysis and
  applications \textbf{20} (1998), no.~2, 513--561.

\bibitem{Benson06}
Hande~Y. Benson, Arun Sen, David~F. Shanno, and Robert~J. Vanderbei,
  \emph{Interior-point algorithms, penalty methods, and equilibrium problems},
  Computational Optimization and Applications \textbf{34} (2006), no.~2,
  155--182.

\bibitem{Bento09}
Antonio~M. Bento, Lawrence~H. Goulder, Mark~R. Jacobsen, and Roger~H. von
  Haefen, \emph{Distributional and efficiency impacts of increased us gasoline
  taxes}, American Economic Review \textbf{99} (2009), no.~3, 667--699.

\bibitem{Benzi05}
Michele Benzi, Gene~H Golub, and Jorg Leisen, \emph{Numerical solution of
  saddle point problems}, Acta Numerica \textbf{14} (2005), 1--137.

\bibitem{Berry95}
Steven Berry, James Levinsohn, and Ariel Pakes, \emph{Automobile prices in
  market equilibrium}, Econometrica \textbf{63} (1995), no.~4, 841--890.

\bibitem{Berry04a}
\bysame, \emph{Differentiated products demand systems from a combination of
  micro and macro data: The new car market}, Journal of Political Economy
  \textbf{112} (2004), no.~1, 68--105.

\bibitem{Bertsekas99}
Dimitri~P. Bertsekas, \emph{Nonlinear programming}, Athena Scientific, 1999.

\bibitem{Besanko98}
David Besanko, Sachin Gupta, and Dipak Jain, \emph{Logit demand estimation
  under competitive pricing behavior: An equilibrium framework}, Management
  Science \textbf{44} (1998), no.~11, 1533--1547.

\bibitem{BLAS}
BLAS, \emph{Basic linear algebra subprograms}, 2010.

\bibitem{Brown90}
Peter~N. Brown and Youcef Saad, \emph{Hybrid krylov methods for nonlinear
  systems of equations}, SIAM Journal of Scientific and Statistical Computing
  \textbf{11} (1990), no.~3, 450--481.

\bibitem{Bunch77}
James~R. Bunch and Linda Kaufman, \emph{Some stable methods for calculating
  inertia and solving symmetric linear systems}, Mathematics of Computation
  \textbf{31} (1977), no.~137, 163--179.

\bibitem{Byrd00}
Richard~H. Byrd, Jean~Charles Gilbert, and Jorge Nocedal, \emph{A trust region
  method based on interior point techniques for nonlinear programming},
  Mathematical Programming, Series A \textbf{89} (2000), no.~1, 149--185.

\bibitem{Byrd04}
Richard~H. Byrd, Nicholas I.~M. Gould, Jorge Nocedal, and Richard~A. Waltz,
  \emph{An algorithm for nonlinear optimization using linear programming and
  quality constrained subproblems}, Mathematical Programming, Series B
  \textbf{100} (2004), no.~1, 27--48.

\bibitem{Byrd99}
Richard~H. Byrd, Mary~E. Hribar, and Jorge Nocedal, \emph{An interior point
  algorithm for large-scale nonlinear programming}, SIAM Journal on
  Optimization \textbf{9} (1999), no.~4, 877--900.

\bibitem{Byrd06}
Richard~H. Byrd, Jorge Nocedal, and Richard~A. Waltz, \emph{{KNITRO}: An
  integrated package for nonlinear optimization}, Tech. report, Ziena Inc.,
  2006.

\bibitem{LANCELOT}
Andrew~R. Conn, Nicholas I.~M. Gould, and Philippe~L. Toint, \emph{Lancelot}.

\bibitem{Conn92}
\bysame, \emph{{LANCELOT}: A fortran package for large-scale nonlinear
  optimization (release a).}, Springer Series in Computational Mathematics,
  vol.~17, Springer-Verlag, 1992.

\bibitem{Conn00}
\bysame, \emph{Trust region methods}, SIAM, 2000.

\bibitem{LA15}
Science \& Technology~Facilities Council, \emph{{HSL} mathematical software
  library catalogue: {LA}15 v 1.2.0}, Tech. report, Research Councils UK, 2010.

\bibitem{Dennis96}
John~E. Dennis and Robert~B. Schnabel, \emph{Numerical methods for
  unconstrained optimization and nonlinear equations}, SIAM, 1996.

\bibitem{Dirkse95}
Steven~P. Dirkse and Michael~C. Ferris, \emph{The {PATH} solver: A non-monotone
  stabilization scheme for mixed complementarity problems}, Optimization
  Methods and Software \textbf{5} (1995), 123--156.

\bibitem{Dolan04}
Elizabeth~D. Dolan, Jorge~J. More, and Todd~S. Munson, \emph{Benchmarking
  optimization software with {COPS} 3.0}, Tech. Report ANL/MCS-273, Argonne
  National Laboratory, February 2004.

\bibitem{Duff04}
Iain~S. Duff, \emph{{MA57} - a code for the solution of sparse symmetric
  definite and indefinite systems}, ACM Transactions on Mathematical Software
  \textbf{30} (2004), no.~2, 118--144.

\bibitem{Ehrenmann09}
Andreas Ehrenmann and Karsten Neuhoff, \emph{A comparison of electricity market
  designs in networks}, Operations Research \textbf{57} (2009), no.~2,
  274--286.

\bibitem{Eisenstat98}
Stanley~C. Eisenstat and Ilse C.~F. Ipsen, \emph{Relative perturbation results
  for eigenvalues and eigenvectors of diagonalisable matrices}, SIAM Journal of
  Matrix Analysis and Applications \textbf{20} (1998), no.~1, 149--158.

\bibitem{Facchinei09}
Francisco Facchinei and Christian Kanzow, \emph{Generalized nash equilibrium
  problems}, Annals of Operations Research \textbf{175} (2009), 177--211.

\bibitem{Ferris97}
Michael~C. Ferris and Jong-Shi Pang, \emph{Engineering and economic
  applications of complementarity problems}, SIAM Review \textbf{39} (1997),
  no.~4, 669--713.

\bibitem{Fletcher87}
Roger Fletcher, \emph{Practical methods of optimization}, Wiley and Sons, 1987.

\bibitem{Fletcher02b}
Roger Fletcher, Nicholas I.~M. Gould, Sven Leyffer, Philippe~L. Toint, and
  Andreas Wachter, \emph{Global convergence of a trust-region sqp-filter
  algorithm for general nonlinear programming}, SIAM Journal on Optimization
  \textbf{13} (2002), no.~2, 635--659.

\bibitem{Fletcher02a}
Roger Fletcher and Sven Leyffer, \emph{Nonlinear programming without a penalty
  function}, Mathematical Programming, Series A \textbf{91} (2002), no.~2,
  239--269.

\bibitem{Fletcher02c}
Roger Fletcher, Sven Leyffer, and Philippe~L. Toint, \emph{On the global
  convergence of a filter-{SQP} algorithm}, SIAM Journal on Optimization
  \textbf{13} (2002), no.~1, 44--59.

\bibitem{Forsgren07}
Anders Forsgren, Philip~E. Gill, and Joshua~D. Griffin, \emph{Iterative
  solution of augmented systems arising in interior methods}, SIAM Journal on
  Optimization \textbf{18} (2007), no.~2, 666--690.

\bibitem{Frischknecht10}
Bart~D. Frischknecht, Katie Whitefoot, and Panos~Y. Papalambros, \emph{On the
  {S}uitability of {E}conometric {D}emand {M}odels in {D}esign for {M}arket
  {S}ystems}, ASME Journal of Mechanical Design \textbf{132} (2010),
  no.~121007, 1--11.

\bibitem{Gabriel05}
Steven~A. Gabriel, Supat Kiet, and Jifang Zhuang, \emph{A mixed
  complementarity-based equilibrium model of natural gas markets}, Operations
  Research \textbf{53} (2005), no.~5, 799--818.

\bibitem{Gabriel01}
Steven~A. Gabriel, Andy~S. Kydes, and Peter Whitman, \emph{The national energy
  modeling system: A large-scale energy-economic equilibrium model}, Operations
  Research \textbf{49} (2001), no.~1, 14--25.

\bibitem{Gill87}
P.~E. Gill, W.~Murray, M.~A. Saunders, and M.~H. Wright, \emph{Maintaining {LU}
  factors of a general sparse matrix}, Linear Algebra and its Applications
  \textbf{88/89} (1987), 239--270.

\bibitem{Gill05}
Philip~E. Gill, Walter Murray, and Michael~A. Saunders, \emph{{SNOPT}: An {SQP}
  algorithm for large-scale constrained optimization}, SIAM Review \textbf{47}
  (2005), no.~1, 99--131.

\bibitem{Gill91}
Philip~E. Gill, Walter Murray, Michael~A. Saunders, and Margaret~H. Wright,
  \emph{Inertia-controlling methods for general quadratic programming}, SIAM
  Review \textbf{33} (1991), no.~1, 1--36.

\bibitem{Gill01}
\bysame, \emph{User's guide for {NPSOL} 5.0: A {FORTRAN} package for nonlinear
  programming}, Tech. Report SOL 86-6, Stanford University, 2001.

\bibitem{Gill81}
Philip~E. Gill, Walter Murray, and Margaret~H. Wright, \emph{Practical
  optimization}, Academic Press, 1981.

\bibitem{Gillespie51}
Robert~Pollack Gillespie, \emph{Partial differentiation}, Interscience
  Publishers, 1951.

\bibitem{Goldberg95}
Pinelopi~K. Goldberg, \emph{Product differentiation and oligopoly in
  international markets: The case of the u.s. automobile industry},
  Econometrica \textbf{63} (1995), no.~4, 891--951.

\bibitem{Goldberg98}
\bysame, \emph{The effects of the corporate average fuel efficiency standards
  in the {US}}, The Journal of Industrial Economics \textbf{46} (1998), no.~1,
  1--33.

\bibitem{Golub96}
Gene~H Golub and Charles F~Van Loan, \emph{Matrix computations}, The Johns
  Hopkins University Press, 1996.

\bibitem{Goolsbee04}
Austan Goolsbee and Amil Petrin, \emph{The consumer gains from direct broadcast
  satellites and the competition with cable tv}, Econometrica \textbf{72}
  (2004), no.~2, 351--381.

\bibitem{Gould85}
Nicholas I.~M. Gould, \emph{On practical conditions for the existence and
  uniqueness of solutions to the general equality quadratic programming
  problems}, Mathematical Programming \textbf{32} (1985), 90--99.

\bibitem{Gould01}
Nicholas I.~M. Gould, Mary~E. Hribar, and Jorge Nocedal, \emph{On the solution
  of equality constrained quadratic programming problems arising in
  optimization}, SIAM Journal of Scientific Computing \textbf{23} (2001),
  no.~4, 1376--1395.

\bibitem{Han77}
S.~P. Han, \emph{A globally convergent method for nonlinear programming},
  Journal of Optimization Theory and Applications \textbf{22} (1977), no.~3,
  297--309.

\bibitem{Higham02}
Nicholas~J. Higham, \emph{Accuracy and stability of numerical algorithms}, 2nd
  edition ed., SIAM, 2002.

\bibitem{Hobbs07}
Benjamin~F. Hobbs and J.~S. Pang, \emph{Nash-cournot equilibria in electric
  power markets with piecewise linear demand functions and joint constraints},
  Operations Research \textbf{55} (2007), no.~1, 113--127.

\bibitem{Hu07}
Xinmin Hu and Daniel Ralph, \emph{Using {EPEC}s to model bilevel games in
  restructured electricity markets with locational prices}, Operations Research
  \textbf{55} (2007), no.~5, 809--827.

\bibitem{Hu07b}
Xinmin Hu, Daniel Ralph, Eric~K. Ralph, Peter Bardsley, and Michael~C. Ferris,
  \emph{Electricity generation with looped transmission networks: Bidding to an
  iso}, 2007.

\bibitem{Jacobsen10}
Mark~R. Jacobsen, \emph{Evaluating u. s. fuel economy standards in a model with
  producer and household heterogeneity}, Working Paper, Stanford University,
  January 2010.

\bibitem{Judd98}
Kenneth~L. Judd, \emph{Numerical methods in economics}, MIT Press, 1998.

\bibitem{Kleit04}
Andrew Kleit, \emph{Impacts of long-range increases in the corporate average
  fuel economy (cafe) standard}, Economic Inquiry \textbf{42} (2004), no.~2,
  279--294.

\bibitem{LAPACK}
LAPACK, \emph{Lapack - linear algebra package}, 2010.

\bibitem{Luenberger09}
David~G. Luenberger and Yinyu Ye, \emph{Linear and nonlinear programming},
  third edition ed., Springer, 2009.

\bibitem{McFadden78}
Daniel~L. McFadden, \emph{Appendix {A}: Definite quadratic forms subject to
  constraint}, Production Economics: A Dual Approach to Theory and Applications
  (Melvyn Fuss and Daniel~L. McFadden, eds.), vol. I: The Theory of Production,
  Amsterdam: North-Holland, 1978.

\bibitem{Michalek04}
Jeremy~J. Michalek, Panos~Y. Papalambros, and Steven~J. Skerlos, \emph{A
  {S}tudy of {F}uel {E}fficiency and {E}mission {P}olicy {I}mpact on {O}ptimal
  {V}ehicle {D}esign {D}ecisions}, ASME Journal of Mechanical Design
  \textbf{126} (2004), 1062--1070.

\bibitem{Morrow10}
W.~Ross Morrow and Steven~J. Skerlos, \emph{Fixed-point approaches to computing
  bertrand-nash equilibrium prices under mixed-logit demand}, Operations
  Research (Forthcoming).

\bibitem{Munson00}
Todd~S. Munson, \emph{Algorithms and environments for complementarity}, Ph.D.
  thesis, University of Wisconsin, Madison, 2000.

\bibitem{Murtagh78}
Bruce~A. Murtagh and Michael~A. Saunders, \emph{Large-scale linearly
  constrained optimization}, Mathematical Programming \textbf{14} (1978),
  no.~1, 41--72.

\bibitem{Murtagh82}
\bysame, \emph{A projected lagrangian algorithm and its implementation for
  sparse nonlinear constraints}, Mathematical Programming Study \textbf{16}
  (1982), 84--117.

\bibitem{Murtagh98}
\bysame, \emph{{MINOS} 5.5 user's guide}, Tech. Report SOL 83-20R, Stanford
  University, 1998.

\bibitem{Nash98}
Stephen~G. Nash, \emph{{N}onlinear {P}rogramming}, ORMS Today (1998).

\bibitem{Nevo00a}
Aviv Nevo, \emph{Mergers with differentiated products: The case of the
  ready-to-eat cereal industry}, The RAND Journal of Economics \textbf{31}
  (2000), no.~3, 395--421.

\bibitem{Nevo01}
\bysame, \emph{Measuring market power in the ready-to-eat cereal industry},
  Econometrica \textbf{69} (2001), no.~2, 307--342.

\bibitem{Nocedal06}
Jorge Nocedal and Stephen~J. Wright, \emph{Numerical optimization},
  Springer-Verlag, 2006.

\bibitem{Papalambros00}
Panos~Y. Papalambros and Douglass~J. Wilde, \emph{Principles of optimal design:
  Modeling and computation}, Cambridge University Press, 2000.

\bibitem{Pernice98}
Michael Pernice and Homer~F. Walker, \emph{{NITSOL}: A newton iterative solver
  for nonlinear systems}, SIAM Journal of Scientific Computing \textbf{19}
  (1998), no.~1, 302--318.

\bibitem{Petrin02}
Amil Petrin, \emph{Quantifying the benefits of new products: The case of the
  minivan}, Journal of Political Economy \textbf{110} (2002), no.~4, 705--729.

\bibitem{Powell78b}
M.~J.~D. Powell, \emph{The convergence of variable matric methods for
  nonlinearly constrained optimization calculations}, Nonlinear Programming
  (O.~L. Mangasarian, R.R. Meyer, and S.M. Robinson, eds.), vol.~3, Academic
  Press, 1978.

\bibitem{Powell78a}
\bysame, \emph{A fast algorithm for nonlinearly constrained optimization
  calculations}, Numerical Analysis (G.~A. Watson, ed.), Lecture Notes in
  Mathematics, vol. 630, Springer-Verlag, 1978.

\bibitem{Ralph94}
Daniel Ralph, \emph{Global convergence of damped newton's method for nonsmooth
  equations, via the path search}, Mathematics of Operations Research
  \textbf{19} (1994), 352--389.

\bibitem{LUSOL}
Michael~A. Saunders, \emph{{LUSOL}: Sparse {LU} for {A}x = b}.

\bibitem{Schenk04}
Olaf Schenk and Klaus Gartner, \emph{Solving unsymmetric sparse systems of
  linear equations with {PARDISO}}, Journal of Future Generation Computer
  Systems \textbf{20} (2004), no.~3, 475--487.

\bibitem{Schenk06}
\bysame, \emph{On fast factorization pivoting methods for symmetric indefinite
  systems}, Electronic Transactions on Numerical Analysis \textbf{23} (2006),
  158--179.

\bibitem{Schenk08}
Olaf Schenk, Andreas Wachter, and Martin Weiser, \emph{Inertia-revealing
  preconditioning for large-scale nonconvex constrained optimization}, SIAM
  Journal of Scientific Computing \textbf{31} (2008), no.~2, 939--960.

\bibitem{Shiau09a}
Ching-Shin~Norman Shiau and Jeremy~J. Michalek, \emph{Optimal {P}roduct
  {D}esign {U}nder {P}rice {C}ompetition}, ASME Journal of Mechanical Design
  \textbf{131} (2009), no.~071003, 1--10.

\bibitem{Shiau09b}
\bysame, \emph{Should {D}esigners {W}orry {A}bout {M}arket {S}tructure?}, ASME
  Journal of Mechanical Design \textbf{131} (2009), no.~011011, 1--9.

\bibitem{Shiau09c}
Ching-Shin~Norman Shiau, Jeremy~J. Michalek, and Chris~T. Hendrickson, \emph{A
  {S}tructural {A}nalysis of {V}ehicle {D}esign {R}esponses to {C}orporate
  {A}verage {F}uel {E}conomy {P}olicy}, Transportation Research A: Policy and
  Practice \textbf{43} (2009), 814--828.

\bibitem{Skerlos05}
Steven~J. Skerlos, Jeremy~J. Michalek, and W.~Ross Morrow, \emph{Sustainable
  {D}esign {E}ngineering and {S}cience: {S}elected {C}hallenges and {C}ase
  {S}tudies}, Sustainability Science and Engineering, Volume 1: Defining
  Principles (M.~A. Abraham, ed.), Elsevier Science, 2005.

\bibitem{Smith04}
Howard Smith, \emph{Supermarket choice and supermarket competition in market
  equilibrium}, The Review of Economic Studies \textbf{71} (2004), 235--263.

\bibitem{Sudhir01}
K~Sudhir, \emph{Competitive pricing behavior in the auto market: A structural
  analysis}, Marketing Science \textbf{20} (2001), no.~1, 42--60.

\bibitem{Thomadsen05}
Raphael Thomadsen, \emph{The effect of ownership structure on prices in
  geographically differentiated markets}, The RAND Journal of Economics
  \textbf{36} (2005), no.~4, 908--929.

\bibitem{Train03}
Kenneth Train, \emph{Discrete choice methods with simulation}, Cambridge
  University Press, 2003.

\bibitem{Trefethen97}
Lloyd~N. Trefethen and David Bau, \emph{Numerical linear algebra}, SIAM, 1997.

\bibitem{Vanderbei06}
Robert~J. Vanderbei, \emph{Loqo user's manual - version 4.05}, Tech. Report
  ORFE-99-??, Princeton University, 2006.

\bibitem{Vanderbei99}
Robert~J. Vanderbei and David~F. Shanno, \emph{An interior point algorithm for
  nonconvex nonlinear programming}, Computational Optimization and Applications
  \textbf{13} (1999), no.~1-3, 231--252.

\bibitem{Wachter06}
Andreas Wachter and Lorenz~T. Biegler, \emph{On the implementation of an
  interior-point filter line search algorithm for large-scale nonlinear
  programming}, Mathematical Programming \textbf{106} (2006), no.~1, 25--57.

\end{thebibliography}
